\documentclass{article}
\usepackage[english]{babel}
\usepackage[utf8]{inputenc}
\usepackage{amsmath}
\usepackage{graphicx}
\usepackage{amsfonts}
\usepackage{enumitem}
\usepackage{amssymb}
\usepackage[colorinlistoftodos]{todonotes}
\usepackage{geometry}
\usepackage{amsthm}
\usepackage{comment}
\usepackage{hyperref}
\usepackage[initials,nobysame]{amsrefs}
\usepackage{cleveref}
\usepackage{tikz}
\usetikzlibrary{math, calc}
\usepackage[final]{showlabels}

\newcommand{\abs}[1]{\left\lvert#1\right\rvert}

\newcommand{\norm}[1]{\left\lVert#1\right\rVert}

\newcommand{\ang}[1]{\left\langle #1 \right\rangle}
\newcommand{\floor}[1]{\left\lfloor #1 \right\rfloor}
\newcommand{\ceil}[1]{\left\lceil #1 \right\rceil}
\newcommand{\paren}[1]{\left( #1 \right)}
\newcommand{\sqb}[1]{\left[ #1 \right]}
\newcommand{\set}[1]{\left\{ #1 \right\}}
\newcommand{\setcond}[2]{\left\{ #1 \;\middle\vert\; #2 \right\}}

\newcommand{\one}{\mathbb{I}}

\newcommand{\ol}{\overline}

\newcommand{\eps}{\epsilon}
\newcommand{\dD}{\delta}

\newcommand{\EE}{\mathbb{E}}

\newcommand{\RR}{\mathbb{R}}
\newcommand{\PP}{\mathbb{P}}
\newcommand{\NN}{\mathbb{N}}
\newcommand{\ZZ}{\mathbb{Z}}

\newcommand{\cE}{\mathcal{E}}
\newcommand{\cF}{\mathcal{F}}

\newcommand{\cH}{\mathcal{H}}

\newcommand{\cL}{\mathcal{L}}
\newcommand{\cM}{\mathcal{M}}
\newcommand{\cN}{\mathcal{N}}

\newcommand{\cP}{\mathcal{P}}
\newcommand{\cQ}{\mathcal{Q}}

\newcommand{\cS}{\mathcal{S}}
\newcommand{\cT}{\mathcal{T}}

\newcommand{\cW}{\mathcal{W}}

\DeclareMathOperator{\conv}{conv}

\newtheorem{thm}{Theorem}[section]
\newtheorem{lem}[thm]{Lemma}

\newtheorem{conj}[thm]{Conjecture}

\newtheorem{prob}[thm]{Problem}
\newtheorem*{rmk}{Remark}

\theoremstyle{definition}
\newtheorem{defn}[thm]{Definition}

\crefname{thm}{theorem}{theorems}
\Crefname{thm}{Theorem}{Theorems}
\crefname{lem}{lemma}{lemmas}
\Crefname{lem}{Lemma}{Lemmas}
\crefname{cor}{corollary}{corollaries}
\Crefname{cor}{Corollary}{Corollaries}
\crefname{conj}{conjecture}{conjectures}
\Crefname{conj}{Conjecture}{Conjectures}
\crefname{claim}{claim}{claims}
\Crefname{claim}{Claim}{Claims}
\crefname{obs}{observation}{observations}
\Crefname{obs}{Observation}{Observations}
\crefname{prop}{proposition}{propositions}
\Crefname{prop}{Proposition}{Propositions}
\crefname{prob}{problem}{problems}
\Crefname{prob}{Problem}{Problems}
\crefname{defn}{definition}{definitions}
\Crefname{defn}{Definition}{Definitions}
\crefname{example}{example}{examples}
\Crefname{example}{Example}{Examples}

\AddToHook{env/thm/begin}{\crefalias{thm}{thm}}
\AddToHook{env/lem/begin}{\crefalias{thm}{lem}}
\AddToHook{env/cor/begin}{\crefalias{thm}{cor}}
\AddToHook{env/conj/begin}{\crefalias{thm}{conj}}
\AddToHook{env/claim/begin}{\crefalias{thm}{claim}}
\AddToHook{env/obs/begin}{\crefalias{thm}{obs}}
\AddToHook{env/prop/begin}{\crefalias{thm}{prop}}
\AddToHook{env/prob/begin}{\crefalias{thm}{prob}}
\AddToHook{env/defn/begin}{\crefalias{thm}{defn}}
\AddToHook{env/example/begin}{\crefalias{thm}{example}}

\DeclareMathOperator{\cost}{cost}
\DeclareMathOperator{\supp}{supp}
\DeclareMathOperator{\msupp}{msupp}

\title{Dominoes tiling a common figure}
\author{Peter Keevash\thanks{Mathematical Institute, University of Oxford, UK. Email: {\tt keevash@maths.ox.ac.uk}. Research supported by ERC Advanced Grant 883810.} \and Jeck Lim\thanks{Institute for Applied Mathematics, University of Bonn, 53115 Bonn, Germany. Email: \texttt{jlim@uni-bonn.de}. Supported by the Deutsche Forschungsgemeinschaft (DFG, German Research Foundation) under Germany's Excellence Strategy -- EXC-2047/2 -- 390685813. Part of this work was done while the author was at the University of Oxford supported by ERC Advanced Grant 883810.}}
\date{}

\begin{document}
\maketitle

\begin{abstract}
    Call a two-element subset of $\mathbb{Z}^2$ a domino. Given finitely many dominoes $D_1,\ldots,D_k$, we show that there exists a finite non-empty $B\subset \mathbb{Z}^2$ such that $B$ can be partitioned into disjoint translates of $D_i$ for all $i$. The proof uses a new absorption mechanism for perfect matchings in algebraically defined hypergraphs.
\end{abstract}


\section{Introduction} \label{sec:intro}

Call a two-element subset of $\ZZ^2$ a \emph{domino}, which does not have to be connected.
We say that a set $D$ \emph{tiles} a set $B$, or that $B$ is \emph{tileable} by $D$,
if $B$ can be partitioned into disjoint translates of $D$.
Given a finite collection of dominoes, Sagdeev~\cite{MO} asked whether there is a common finite figure that is simultaneously tileable by all of them.

\begin{prob} \label{prob:dominoes}
    Given dominoes $D_1,\ldots,D_k\subset\ZZ^2$, is there a finite
    non-empty set $B\subset\ZZ^2$ that is tileable by every $D_i$?
\end{prob}

For example, for the two standard dominoes:
\begin{center}
    \begin{tikzpicture}[x=0.45cm, y=0.45cm, line width=0.6pt]
        \draw (0,0) rectangle (2,1);
        \draw (1,0) -- (1,1);
        \node at (1,-0.5) {$D_1$};
        \begin{scope}[shift={(6,0)}]
            \draw (0,0) rectangle (1,2);
            \draw (0,1) -- (1,1);
            \node at (0.5,-0.5) {$D_2$};
        \end{scope}
    \end{tikzpicture}
\end{center}
we can take $B$ to be the $2\times 2$ square. If we also include the two diagonal dominoes:
\begin{center}
    \begin{tikzpicture}[x=0.45cm, y=0.45cm, line width=0.6pt]
        \draw (0,0) rectangle (1,1);
        \draw (1,1) rectangle (2,2);
        \node at (1,-0.5) {$D_3$};
        \begin{scope}[shift={(6,0)}]
            \draw (0,1) rectangle (1,2);
            \draw (1,0) rectangle (2,1);
            \node at (1,-0.5) {$D_4$};
        \end{scope}
    \end{tikzpicture}
\end{center}
we can take $B$ to be the following figure, found by Sagdeev~\cite{MO}.
\begin{center}
    \begin{tikzpicture}[x=0.35cm, y=0.35cm, line width=0.5pt]
        \foreach \y/\xfirst/\xlast in {
                7/2/3, 6/1/4, 5/1/4, 4/0/5, 3/0/5, 2/1/4, 1/1/4, 0/2/3
            } {
                \foreach \x in {\xfirst,...,\xlast}
                \draw (\x,\y) rectangle ++(1,1);
            }
    \end{tikzpicture}
\end{center}
If we also include the horizontal domino with a gap:
\begin{center}
    \begin{tikzpicture}[x=0.45cm, y=0.45cm, line width=0.6pt]
        \draw (0,0) rectangle (1,1);
        \draw (2,0) rectangle (3,1);
        \node at (1.5,-0.5) {$D_5$};
    \end{tikzpicture}
\end{center}
constructing $B$ becomes much more difficult. Using a SAT solver, Tikhomirov \cite[Comments]{MO} found the following solution with 416 cells, which is minimal in the sense that there are no solutions with a smaller bounding box.
\begin{center}
    \begin{tikzpicture}[x=0.24cm, y=0.24cm, line width=0.35pt]
        \foreach \y/\xfirst/\xlast in {
                13/30/33,
                12/15/26,
                12/29/36,
                11/10/37,
                10/6/17,
                10/24/39,
                9/3/42,
                8/1/44,
                7/0/31,
                7/34/45,
                6/0/31,
                6/34/45,
                5/1/44,
                4/3/42,
                3/6/17,
                3/24/39,
                2/10/37,
                1/15/26,
                1/29/36,
                0/30/33
            } {
                \foreach \x in {\xfirst,...,\xlast}
                \draw (\x,\y) rectangle ++(1,1);
            }
    \end{tikzpicture}
\end{center}

We answer Sagdeev's question affirmatively.
In fact, we prove the result more generally for arithmetic progressions of arbitrary lengths.
An \emph{arithmetic progression} (AP) in $\ZZ^2$ is a set of the form
\[
    \{a,a+v,\ldots,a+(d-1)v\},
\]
where $a\in\ZZ^2$, $v\in\ZZ^2\setminus\{0\}$ and $d\geq 2$ is its length.
In particular, dominoes are APs of length two. Our main result is as follows.
\begin{thm} \label{thm:main-ap}
    Let $D_1,\ldots,D_k\subset\ZZ^2$ be APs. Then there exists a
    finite non-empty $B\subset\ZZ^2$ that is tileable by every $D_i$.
\end{thm}

Clearly \Cref{thm:main-ap} implies a positive answer to \Cref{prob:dominoes}.
Sagdeev~\cite[Answer]{MO} showed that in fact they are equivalent: for each AP there is a finite family of dominoes in its direction such that every finite set tiled by each of those dominoes is also tiled by the AP.
We do not use this equivalence in our proof. Instead, we work with APs directly as they are more convenient. More details, including our proof sketch, are given in the next section.

\subsection*{Acknowledgements}
The second author thanks David Conlon for suggesting this problem and for helpful discussions.

\subsection*{Statement on AI use}
An initial draft of the paper was written by the authors by February 2026 without significant AI assistance. Over the course of several months, various versions of GPT (5.4 Plus, 5.5 Plus, 6 Astra, 6.1 Sol) were used for figure generation, proofreading and editing, with no changes to the main arguments of the proof.

\section{Overview of the proof}  \label{sec:overview}

The main idea in our construction is a new absorption mechanism for perfect matchings in algebraically defined hypergraphs.
To explain where this comes into play, we must first set up the overall strategy for the construction.
We start with some preliminary observations.
The generalisation to APs has two benefits. First, we may assume that $D_1,\ldots,D_k$ have different directions. Indeed, if, say, $D_1,D_2$ have the same direction, then there is an AP $D$ in the same direction such that if $B$ is tileable by $D$, then it is tileable by $D_1$ and $D_2$. So we may replace $D_1,D_2$ by the single AP $D$. Second, we may also assume that the direction $v$ of $D_i$ is primitive, i.e., the coordinates of $v$ are coprime. Indeed, if $v=cv'$ for some $c>1$, then we may replace $D_i$ with the AP $D_i'$ of length $cd$ and direction $v'$. We may also assume that all the APs $D_1,\ldots,D_k$ have the same length, by replacing their lengths with their lowest common multiple.

Henceforth, let $v_1,\ldots,v_k$ be the primitive, pairwise non-parallel directions of $D_1,\ldots,D_k$, respectively, and assume that all the APs have length $d$.

For a line $\ell$ parallel to some AP $D_i$, we say that $B$ is \emph{tileable along $\ell$} if $B\cap \ell$ is tileable by $D_i$.
Our aim is to find a finite non-empty set $B$ which is tileable along every line $\ell$ parallel to some $D_i$.

We will construct such a finite set $B$ via the following steps.
\begin{enumerate}
    \item \textbf{Large polygon.} We will obtain $B$ from a series of small modifications to a large symmetric polygon with $2k$ sides parallel to the directions $v_1,\ldots,v_k$. Let $B_1$ be the set of integer points in the polygon $n\cdot \conv(\set{-c_1 v_1, c_1 v_1}+\cdots+\set{-c_k v_k, c_k v_k})$ for some large integer $n$ and some fixed positive integers $c_1,\ldots,c_k$.

    \item \textbf{Random jagged boundary.} We obtain $B_2$ from $B_1$ by attaching random jagged caps of thickness $O(1)$ along the boundary of $B_1$. These caps are constructed using a Markov chain with the properties that (a) for almost any point $p$ and any direction $v_i$, the distance along $v_i$ from $p$ to the boundary is uniformly distributed mod $d$, and (b) sufficiently separated constant-sized sets of events as in (a) are approximately independent. In particular, about $1/d$ of the lines in directions $v_1,\dots,v_k$ intersecting $B_2$ are lines along which $B_2$ is already tileable. We call these lines \emph{pure}, which will serve as a large reserve for the subsequent cover down, absorption, and gadget steps.

    \item \textbf{Almost-perfect matching.} For each line $\ell$ along which $B_2$ is not tileable, we partition $\ell\cap B_1$ into segments, aiming to poke exactly one hole in each segment (at certain positions modulo $d$) in order to make $B_2$ tileable along that line. We build a $k$-uniform $k$-partite hypergraph $G$, whose vertices are these segments, with edges corresponding to cells in $B_2$ where the segment requirements allow us to poke holes. A perfect matching in this auxiliary hypergraph would therefore give a set of cells whose removal will give the desired tileable set (a common device in combinatorial construction, see the survey \cite{keevash2019hypergraph}). Such perfect matching problems are generally hard and often approached via the two-step strategy of (i) finding an almost-perfect matching, then (ii) applying an absorption method to improve this to a perfect matching. The main difficulty (in general and in this paper) lies in absorption, although, as explained in \Cref{sec:match}, even the almost-perfect matching step is quite difficult here.

    \item \textbf{Cover down.} The almost-perfect matching step pokes holes in $B_2$ to produce $B_3$ which is tileable along all but a small number of lines.
    To prepare for absorption, we obtain $B_4$ by poking further holes with a more convenient structure for these problematic lines:
    we arrange that they are all pure and pass near the centre of our shape.

    \item \textbf{Tile divisibility.} The main difficulty in absorption is modifying $B_4$ to obtain $B_5$
    that satisfies the integral relaxation of tileability, which we call \emph{tile divisibility}.
    To understand this for a general hypergraph, one can think of the perfect matching problem as finding a $\{0,1\}$-solution $x$ to $Ax=1$,
    where $A$ is the vertex/edge incidence matrix and $1$ is the all-1 vector indexed by the vertices;
    then the integral relaxation is an easier problem that only asks for an integral solution $x$.
    In our context, divisibility is best understood by identifying any finite $S \subset \ZZ^2$
    with the Laurent polynomial $\sum_{(a_1,a_2) \in S} x_1^{a_1} x_2^{a_2} \in \ZZ[x_1^\pm,x_2^\pm]=\ZZ[x_1,x_1^{-1},x_2,x_2^{-1}]$.
    We say that $S$ is divisible along a line $\ell$ in direction $v_i$ if $S \cap \ell$ is divisible by $D_i$ under this identification.
    We say that $S$ is tile divisible if it is divisible along all lines in the directions $D_1,\dots,D_k$.
    This is equivalent to $S$ being divisible by $D_1,\dots,D_k$.
    Applying our new absorption method, discussed below, we obtain $B_5$ from $B_4$ by poking a small number of holes.
    The resulting set $B_5$ is tile divisible and is tileable along all but a small number of lines.

    \item \textbf{Moving holes.} In this final step, we turn tile-divisible lines into tileable lines. Thus we finish resolving the untileable lines left by the almost-perfect matching step (so we do not need to actually find a perfect matching in the auxiliary hypergraph considered above). The intuition is that if $B_5$ is divisible but not tileable along some line $\ell$, then we could make $B_5$ tileable by certain moves of the holes in steps of multiples of $d$ along $\ell$. A naive implementation would cause trouble along other lines, but by using certain (rather complicated) gadgets, this can be done without destroying tileability elsewhere. This gives us the final set $B_6$, which is tileable.
\end{enumerate}

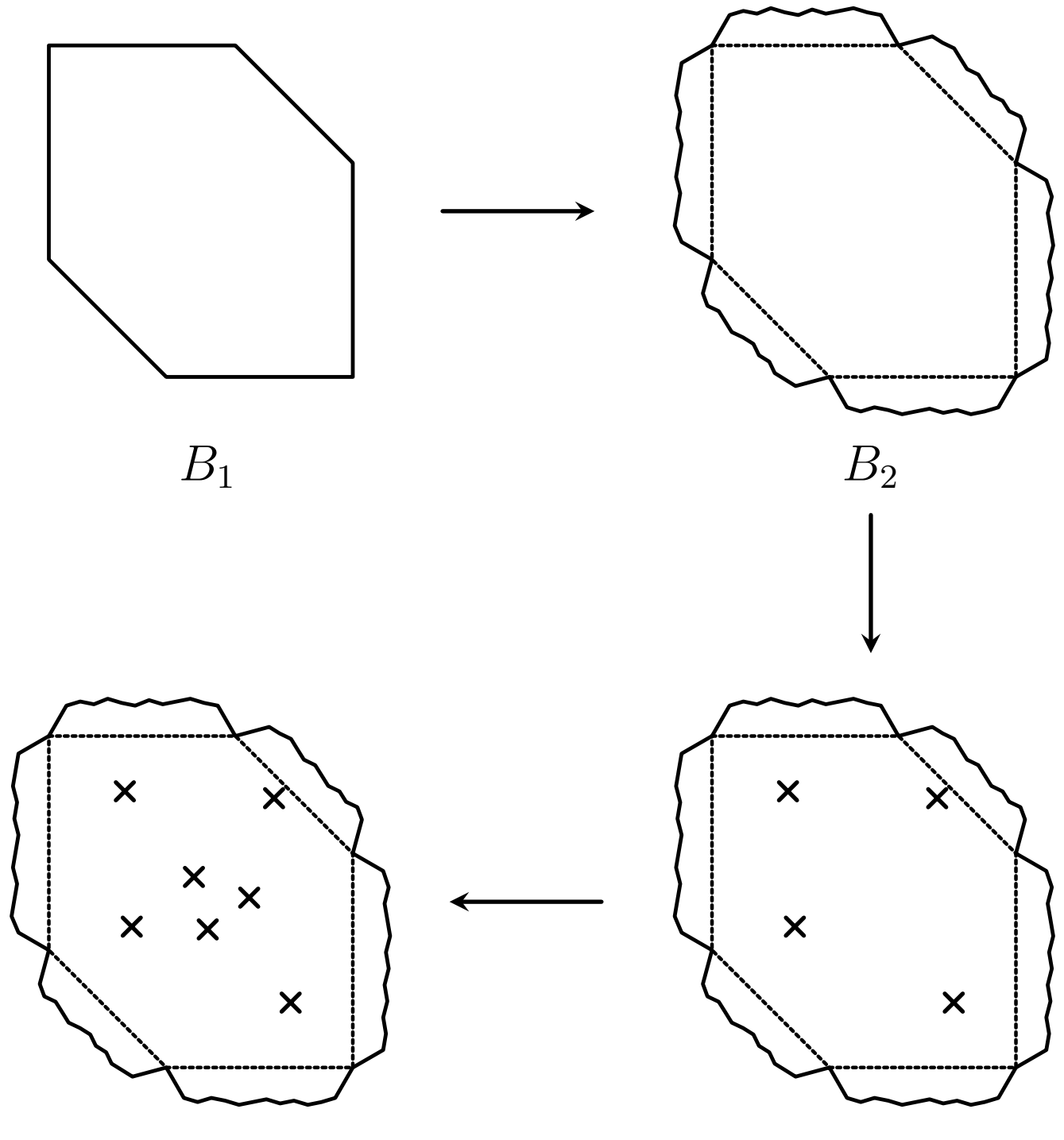
\begin{figure}
    \centering
    \begin{tikzpicture}[
            x=1.15cm, y=1.15cm, line cap=round, line join=round,
            outline/.style={line width=0.9pt},
            original boundary/.style={line width=0.85pt, densely dotted},
            stage arrow/.style={->, >=stealth, line width=1pt},
            stage label/.style={font=\large},
            pics/hole/.style={code={
                            \draw[line width=1pt] (-0.065,-0.065) -- (0.065,0.065)
                            (-0.065,0.065) -- (0.065,-0.065);
                        }}
        ]
        \def\sketchCore{
            (0.10,2.50) -- (1.45,2.50) -- (2.30,1.65)
            -- (2.30,0.10) -- (0.95,0.10) -- (0.10,0.95) -- cycle
        }
        \def\sketchJagged{
            (0.10,2.50)
            -- (0.227,2.720) -- (0.327,2.750) -- (0.426,2.730)
            -- (0.526,2.770) -- (0.626,2.740) -- (0.725,2.720)
            -- (0.825,2.760) -- (0.924,2.730) -- (1.024,2.750)
            -- (1.124,2.770) -- (1.223,2.740) -- (1.323,2.720)
            -- (1.45,2.50)
            -- (1.695,2.566) -- (1.770,2.519) -- (1.853,2.479)
            -- (1.899,2.404) -- (1.946,2.329) -- (2.028,2.289)
            -- (2.075,2.214) -- (2.122,2.139) -- (2.204,2.099)
            -- (2.251,2.024) -- (2.333,1.985) -- (2.366,1.895)
            -- (2.30,1.65)
            -- (2.520,1.523) -- (2.560,1.405) -- (2.530,1.287)
            -- (2.550,1.170) -- (2.570,1.052) -- (2.540,0.934)
            -- (2.560,0.816) -- (2.530,0.698) -- (2.550,0.580)
            -- (2.520,0.463) -- (2.540,0.345) -- (2.520,0.227)
            -- (2.30,0.10)
            -- (2.173,-0.120) -- (2.073,-0.150) -- (1.974,-0.170)
            -- (1.874,-0.140) -- (1.774,-0.160) -- (1.675,-0.130)
            -- (1.575,-0.150) -- (1.476,-0.170) -- (1.376,-0.140)
            -- (1.276,-0.120) -- (1.177,-0.150) -- (1.077,-0.120)
            -- (0.95,0.10)
            -- (0.705,0.034) -- (0.630,0.081) -- (0.554,0.128)
            -- (0.515,0.210) -- (0.440,0.257) -- (0.400,0.339)
            -- (0.325,0.386) -- (0.243,0.425) -- (0.196,0.501)
            -- (0.149,0.576) -- (0.067,0.615) -- (0.034,0.705)
            -- (0.10,0.95)
            -- (-0.120,1.077) -- (-0.170,1.195) -- (-0.150,1.313)
            -- (-0.130,1.430) -- (-0.160,1.548) -- (-0.140,1.666)
            -- (-0.120,1.784) -- (-0.150,1.902) -- (-0.130,2.020)
            -- (-0.160,2.137) -- (-0.140,2.255) -- (-0.120,2.373)
            -- cycle
        }
        \tikzset{
            pics/jagged region/.style={code={
                            \draw[outline] \sketchJagged;
                            \draw[original boundary] \sketchCore;
                        }},
            pics/first holes/.style={code={
                            \foreach \x/\y in {0.65/2.10,1.73/2.05,0.70/1.12,1.85/0.57}
                            \pic at (\x,\y) {hole};
                        }}
        }

        \draw[outline] \sketchCore;
        \node[stage label] at (1.25,-0.55) {$B_1$};

        \begin{scope}[shift={(4.8,0)}]
            \pic {jagged region};
            \node[stage label] at (1.25,-0.55) {$B_2$};
        \end{scope}

        \begin{scope}[shift={(4.8,-5)}]
            \pic {jagged region};
            \pic {first holes};
            \node[stage label] at (1.25,-0.65) {$B_3,\,B_4$};
        \end{scope}

        \begin{scope}[shift={(0,-5)}]
            \pic {jagged region};
            \pic {first holes};
            \foreach \x/\y in {1.15/1.48,1.25/1.10,1.55/1.33}
            \pic at (\x,\y) {hole};
            \node[stage label] at (1.25,-0.65) {$B_5,\,B_6$};
        \end{scope}

        \draw[stage arrow] (2.95,1.30) -- (4.05,1.30);
        \draw[stage arrow] (6.05,-0.90) -- (6.05,-1.90);
        \draw[stage arrow] (4.10,-3.70) -- (3.00,-3.70);
    \end{tikzpicture}
    \caption{The rough structure of the different stages of the construction, for the model case $k=3$. The X's from $B_3$ onwards represent holes. The additional holes poked in $B_5,B_6$ are concentrated around the origin.}
    \label{fig:sketch}
\end{figure}

To help visualise the construction, it would be helpful to consider the model case $k=d=3$, with directions $v_1=(1,0),v_2=(1,-1),v_3=(0,-1)$. See \Cref{fig:sketch} for a cartoon depiction of what the various stages of the construction look like for this model case.

\subsection{Polynomial decompositions} \label{subsec:polydecomp}

The main new ingredient in the proof strategy described above is the absorption mechanism required in Step 5.
Here we will motivate and describe this mechanism, indicating the difficulties in its proof,
which is presented in a companion paper \cite{KL26-1} with further applications.
Recall that we identify any finite $S \subset \ZZ^2$  with the Laurent polynomial
$\sum_{(a_1,a_2) \in S} x_1^{a_1} x_2^{a_2} \in \ZZ[x_1^\pm,x_2^\pm]$.
Our goal in Step 5 is to poke a small number of holes in $B_4$ to obtain $B_5$ that is divisible by $D_1,\dots,D_k$.
We need to ensure that the number of holes added is small
(it will be $O(n^{1-\eps})$ for some fixed $\eps>0$) so that we do not run out of space and we can find
suitable gadgets in the next step (for the purpose of this discussion we will gloss over various further
separation properties of the holes needed to apply the approximate independence property of the random caps).

We will achieve divisibility sequentially: at step $i \in [k]$ we poke holes $F_i$
so that $B_4^i = B_4 - F_1 - \dots - F_i$ is divisible by $D_j$ for all $j \le i$.
Let $H_i$ be the set of holes on lines in direction $i$ that are not tileable in $B_4$.
By the cover down step, each such line has exactly one hole in $H_i$ and $B_4 + H_i$ is $D_i$-divisible.
We write $B_4 + H_i = B_4^{i-1} - F_i^*$, where $F_i^* = -H_i - F_1 - \dots   - F_{i-1}$.
We aim to replace $F_i^*$ by $F_i$ so that
(a) $F_i$ (and so $B_4^i=B_4^{i-1}-F_i$) is $D_j$-divisible for all $j<i$,
and (b) $F_i - F_i^*$ (and so $B_4^i = B_4 + H_i - (F_i - F_i^*)$) is $D_i$-divisible.

For the purpose of this discussion, we pass to a toy version of the problem
where we replace each $D_i = \sum_{j=0}^{d-1} x^{jv_i}$ by $x^{dv_i} - 1 =  (x^{v_i} - 1)D_i$.
The basic building block needed for this toy problem is an iterative procedure
for boosting divisibility of some Laurent polynomial $f \in \ZZ[x_1^\pm,x_2^\pm]$:
if $f$ is divisible by $x^{dv_j} - 1$ for all $j<i$
and there is some small $F^*$ so that $f-F^*$ is divisible by $x^{dv_i} - 1$,
then we want to find some small $F$ so that $f-F$ is divisible by $x^{dv_j} - 1$ for all $j \le i$.
Here we are also relaxing from the $\{0,1\}$-problem of choosing a small number of holes
to the integral problem of choosing a small polynomial;
this relaxation can be addressed with combinatorial gadgets, but the simplification of replacing
$D_i$ by $x^{dv_i} - 1$ requires substantial algebraic arguments (see \Cref{subsec:twovar}).

The basic building block can be made very concrete, and its small cases are deceptively simple.
The first non-trivial case is $i=2$ (for $i=1$ we can set $F=F^*$). It is convenient to change variables
from $(x_1,x_2)$ to $(X_1,X_2)$ where $X_1=v_2$ and we choose $X_2$ so that $\det(X_1,X_2)=1$.
Write $v_1 = u_1 X_1 + u_2 X_2$ and $X^u=X_1^{u_1} X_2^{u_2}$.
By assumption $f-F^*$ is divisible by $X_1^d-1$
and $f$ is divisible by $X^{du}-1$.
Writing $\pi$ for reduction modulo $X_1^d-1$,
we have $\pi(F^*)=\pi(f)$ divisible by $\pi(X^{du}-1)=X_2^{du_2}-1$.
Writing $\pi(F^*)=\pi(f)=\sum_{j=0}^{d-1} X_1^j h_j(X_2)$,
we see that each $h_j(X_2)$ is divisible by $X_2^{du_2}-1$.
By repeatedly extracting $\pm 1$ pairs of coefficients,
we can thus write each $h_j(X_2)$ as a small sum of basic polynomials
of the form $X_2^c (1-X_2^{bdu_2})$ with $b,c \in \ZZ$.
For each such summand of $h_j$, include $X_1^j X_2^c(1-X^{bdu})$ in $F$.
Then $\pi(F)=\sum_{j=0}^{d-1}X_1^j h_j(X_2)=\pi(F^*)=\pi(f)$, so $f-F$ is divisible by $X_1^d-1$.
Also $F$ is divisible by $X^{du}-1$, and hence so is $f-F$, as required.

For general $i$, changing variables again to $(X_1,X_2)$ where $X_1=v_i$,
we have $f-F^*$ divisible by $X_1^d-1$ and $f$ divisible by $X^{du^{(j)}}-1$
where $v_j = u^{(j)}_1 X_1 + u^{(j)}_2 X_2$ for each $j<i$.
Since the $v_j$ are non-parallel, $f$ is divisible by $\prod_{j<i} (X^{du^{(j)}}-1)$.
Thus $\pi(F^*)=\pi(f)=\sum_{\ell=0}^{d-1} X_1^\ell h_\ell(X_2)$
with each $h_\ell(X_2)$ divisible by $\prod_{j<i} (X_2^{duc_2}-1)$.
Now the key input is the following theorem from our companion paper \cite{KL26-1}
showing that $h_\ell(X_2)$ can be written as a small sum of basic polynomials
of the form $X_2^c \prod_{j<i} (1-X_2^{b_j du^{(j)}_2})$ with $b_j,c \in \ZZ$.
We then obtain $F$ by lifting each such summand to
$X_2^c \prod_{j<i} (1-X^{b_j du^{(j)}})$.

\begin{thm}[{\cite[Corollary~4.2]{KL26-1}}] \label{thm:poly-decomp}
    Fix positive integers $m_1,\ldots,m_k$. For any $\varepsilon>0$, there exists $\delta>0$ such that the following holds. Let $F=\sum_i f_i x^i\in \ZZ[x]$ be a polynomial of degree at most $n$ such that $\prod_i (1-x^{m_i})\mid F$ and $\sum |f_i|=O(n^{1-\varepsilon})$. Then $F$ can be written as a sum of $O(n^{1-\delta})$-many polynomials of the form
    \[\pm x^a(1-x^{m_1a_1})\cdots (1-x^{m_ka_k}),\]
    where $a,a_1,\ldots,a_k$ are non-negative integers. Furthermore, each such polynomial has degree at most $n$.
\end{thm}

As discussed in \cite{KL26-1}, the cases $k=1$ and $k=2$ of \Cref{thm:poly-decomp}
are straightforward to prove by a greedy algorithm based on simple swaps,
but for $k \ge 3$ there are lower bound constructions showing that no such greedy algorithm exists.
Our proof therefore uses an iterative process based on alternating cancelling and focusing steps,
where the cancelling step relies on a miraculous coupling of signed distributions
on cross polytopes using Legendre polynomials. It also requires some algebra
and properties of cyclotomic polynomials, but these are mostly isolated from the current paper,
in which the arguments are mainly combinatorial.

\subsection{Organisation of the paper and notation}

We implement Steps 1 and 2 in \Cref{sec:polygon},
Step 3 in \Cref{sec:match},
Step 4 in \Cref{sec:cover-down},
Step 5 in \Cref{sec:div},
and Step 6 in \Cref{sec:tile}.
\Cref{sec:conclude} has some concluding remarks.

Throughout the paper we fix positive integers $k,d$ and primitive non-parallel arithmetic progressions $D_1,\ldots,D_k$.
The case $k=1$ is trivial, so we assume $k\geq 2$. Let $v_i$ be the direction of $D_i$, which is primitive.
We may assume that $D_i=\{0,v_i,\ldots,(d-1)v_i\}$. A \emph{direction} will always be one of the vectors $\pm v_1,\ldots,\pm v_k$.
We always assume that $n$ is a sufficiently large integer. Underlying constants in big-$O$ notation may depend on $k,d,v_1,\ldots,v_k$, but not on $n$.

For a quantity $Q$, we say that $Q$ is \emph{small} if $Q=O(n^{1-\varepsilon})$ for some $\varepsilon>0$.

We say that $Q$ is \emph{tiny} if $Q=O(1)$.

For a polynomial $F=\sum_i a_i x^i\in \ZZ[x]$, we define the \emph{size} of $F$ to be $|F|:=\sum_i |a_i|$.

\section{Polygon with random boundary}  \label{sec:polygon}

This section implements Steps 1 and 2 of the proof strategy. We start by fixing some conventions.
All affine changes of coordinates below preserve the lattice $\ZZ^2$, and change Euclidean distances by bounded factors. All distances are Euclidean distances, but in most cases we are only interested in the distance up to constant factors, so any reasonable distance metric will do. By an affine change of coordinates, we may assume without loss of generality that $v_1=(1,0)$. Furthermore, by flipping signs, we may assume that $v_2,\ldots,v_k$ have negative $y$-coordinates, and we order the vectors $v_1,v_2,\ldots,v_k$ in clockwise order.

Consider the symmetric $2k$-gon $\ol{B}=\conv(\{-c_1v_1,c_1v_1\}+\cdots+\{-c_kv_k,c_kv_k\})$ for some fixed positive integers $c_1,\ldots,c_k$, with sides parallel to each $v_i$. Choose the coefficients so that no vertex lies on any of the lines $\RR v_i$ through the origin $O$, and the perpendicular foot from $O$ to the line containing each side lies in the relative interior of that side. Such a choice exists: the centrally symmetric polygon with sides parallel to the $v_i$ and tangent to a circle centred at $O$ has the latter property and is of the above form with positive real coefficients. A sufficiently small rational perturbation of these coefficients preserves this property and avoids the finitely many linear equations placing a vertex on some $\RR v_i$. Scaling then makes all coefficients positive integers. Thus every ray from $O$ in a direction $\pm v_i$ meets the boundary in the interior of a side.

Let $\ol{B}_1= n\cdot \ol{B}$ for a large integer $n$. In particular, every ray from $O$ in a direction $\pm v_i$ meets the boundary at a point that is at distance $\Omega(n)$ from every vertex. Let $A_1,\ldots,A_{2k}\in \ZZ^2$ be the vertices of $\ol{B}_1$ in clockwise order, so that $A_{i+1}-A_i=2nc_i v_i$ (indices modulo $2k$), where $v_{k+i}=-v_i$ and $c_{k+i}=c_i$ for $i=1,\ldots,k$. Set $B_1=\ol{B}_1\cap \ZZ^2$.

We will modify the boundary of $B_1$ by adding a random ``jagged'' cap on each side to obtain $B_2$.

\subsection{Random jagged cap}

We will describe the random cap construction on the side $A_1A_2$; caps on the other sides are constructed similarly and independently. After a suitable affine transformation and renaming of variables, we assume that the side $A_1A_2$ connects the points $(0,0)$ and $(n,0)$, and $B_1$ lies below this line. We denote the jagged cap to be added by $P$, which lies above the $x$-axis and will be constructed in such a way that $P$ is tileable by $D_1$ and is flat enough so that for every line $\ell$ parallel to some $v_i$ ($i\neq 1$), $\ell\cap P$ is a contiguous segment on $\ell$. This ensures that $B_1\cup P$ is not obviously non-tileable by other $D_i$'s near this side.

Let $w,\bar{h}$ be (bounded) positive integers with $\bar{h}\gg w$ and $d\mid w$. The cap $P$ will be made of rectangular blocks of width $w$ and height approximately $\bar{h}$. We assume $n$ is a multiple of $w$. Choose sufficiently large fixed multiples $w_L,w_R$ of $w$, and let $m$ be such that $w_L+w_R+mw=n$.

Let $y_i$ be the $y$-coordinate of $v_i$, for $i=2,\ldots,k$, orienting these vectors so that $y_i<0$. Take independent random variables $Y,Y_2,\ldots,Y_k$, where $Y_i$ is uniform on $\set{0,y_i,2y_i,\ldots,(d-1)y_i}$ and $Y$ is uniform on $[\max_{i=2}^k (-y_i)]$. Let $\cH$ be the distribution of $\bar{h}+Y+Y_2+\cdots+Y_k$. Its support is an interval of positive integers of length $O(1)$ independent of $\bar{h}$, centered near $\bar{h}$.

On the support of $\cH$, consider the lazy Markov chain with transition matrix $K$ given by
\[
    K(h,h')=\frac14\min\set{1,\frac{\cH(h')}{\cH(h)}}\quad\text{if }|h-h'|=1,
    \qquad K(h,h)=1-\sum_{h'\neq h}K(h,h'),
\]
and all other transitions zero. This chain is irreducible, reversible with stationary distribution $\cH$, and aperiodic. Let $h_1,h_2,\ldots,h_m$ be a run of this chain, with initial state $h_1$ having distribution $\cH$. Then $h_j$ also has distribution $\cH$ for all $j$.
\begin{lem} \label{lem:markov-mixing}
    There exists a constant $c=c_{\ref{lem:markov-mixing}}>0$, independent of $\bar{h}$, such that the following holds. Let $J\subseteq [m]$ consist of one or two consecutive indices, and let every index in $I\subseteq [m]$ be at distance at least $t\geq 1$ from $J$. If $F$ is determined by $(h_j)_{j\in J}$ and $E$ is determined by $(h_i)_{i\in I}$, with $\PP(E)>0$, then
    \[
        \PP(F\mid E)=\PP(F)\pm O(e^{-ct}).
    \]
    In particular, if $J=\{j\}$, then for every $h$ in the support of $\cH$,
    \[
        \PP(h_j=h\mid E)=\cH(h) \pm O(e^{-ct}).
    \]
\end{lem}

The proof makes use of the following exponential convergence of irreducible aperiodic Markov chains.

\begin{thm}[{\cite[Theorem 4.9]{LP17}}] \label{thm:exp-mixing}
    Let $K$ be the transition matrix of an irreducible and aperiodic Markov chain with finite state space $X$ and stationary distribution $\pi$. Then there exist constants $C,c>0$ such that
    \[
        \max_{x\in X}\norm{K^t(x,\cdot)-\pi}_{TV}\leq Ce^{-ct}.
    \]
\end{thm}

\begin{proof}[Proof of \Cref{lem:markov-mixing}]
    Since the minimum positive value of $\cH$ is bounded below independently of $\bar{h}$, by \Cref{thm:exp-mixing}, there exists $c>0$ such that, uniformly over states $x,y$ and integers $r\geq 0$,
    \begin{equation} \label{eq:exp-mixing}
        K^r(x,y)=\cH(y)\bigl(1+O(e^{-cr})\bigr).
    \end{equation}
    Write $J=\{a,\ldots,b\}$, where $b-a\leq 1$, and first condition on all the variables indexed by $I$. By the Markov property, the conditional law on $J$ depends only on the nearest observed index $l<a$ and the nearest observed index $r>b$ of $I$, when these exist. If both exist, then for any states $x,y$ and any tuple $(z_a,\ldots,z_b)$,
    \[
        \PP(h_q=z_q\text{ for all }q\in J\mid h_l=x,h_r=y)
        =
        \frac{K^{a-l}(x,z_a)
            \left(\prod_{q=a}^{b-1}K(z_q,z_{q+1})\right)
            K^{r-b}(z_b,y)}{K^{r-l}(x,y)}.
    \]
    For sufficiently large $t$, \eqref{eq:exp-mixing} shows that this is
    \[
        \cH(z_a)\prod_{q=a}^{b-1}K(z_q,z_{q+1})
        \bigl(1+O(e^{-ct})\bigr).
    \]
    The quantity $\cH(z_a)\prod_{q=a}^{b-1}K(z_q,z_{q+1})$ is just $\PP(h_q=z_q\text{ for all }q\in J)$, as required. A similar argument holds if only one or neither of $l,r$ exists.
\end{proof}

$P$ will consist of $m+2$ parts, labelled $P_L,P_1,P_2,\ldots,P_m,P_R$, in this order from left to right, with widths $w_L,w,\ldots,w,w_R$. Each random piece $P_j$ is a rectangular block of width $w$ and height $h_j$.

The end pieces $P_L,P_R$ are made of rectangular blocks of width $w$ with non-negative integer heights. The block at each corner has height zero, the block next to $P_1$ or $P_m$ has height $h_1$ or $h_m$, respectively, and adjacent block heights differ by at most one. Pick any choice of such heights, which is possible for sufficiently large $w_L,w_R$.

$w$ will be chosen sufficiently large relative to the length of the support of $\cH$, so that \Cref{lem:cap-tau-unif} holds later on.

\subsection{Random properties of the cap}

We say that a point $(x,y)\in \ZZ^2$ \emph{lies below $P$} if either $y\leq 0$ or $(x,y)\in P$.

\begin{defn}
    For a point $p\in \ZZ^2$ and a direction $v\in\{-v_2,\ldots,-v_k\}$, so that $v_y>0$, define the \emph{distance from $p$ to $P$ along $v$}, denoted $\tau_{v,P}(p)$, to be the largest integer $t$ such that $p+tv$ lies below $P$. Define $\pi_{v,P}(p)\in \RR$ to be the $x$-coordinate of the intersection of the line $p+\RR v$ with the line $y=\bar{h}$.
\end{defn}
For these directions, the set of integers $t$ such that $p+tv$ lies below $P$ is of the form $(-\infty,T]\cap\ZZ$. This $T$ is precisely $\tau_{v,P}(p)$.

We first show that the distance from any point $p$ to $P$ along the direction $-v_i$ for each $i=2,\ldots,k$ is uniformly distributed mod $d$, provided that $\pi_{-v_i,P}(p)$ lies in $[0,n]$ and is not too close to either end. In the next two lemmas, the margins in $[O(1),n-O(1)]$ exclude $P_L,P_R$ together with a fixed buffer large enough that the relevant exit occurs in the random part $P_1,\ldots,P_m$.

\begin{lem} \label{lem:cap-tau-unif}
    Let $p\in \ZZ^2$ and $i\in \{2,\ldots,k\}$. Suppose that $\pi_{-v_i,P}(p)\in [O(1),n-O(1)]$. Then $\tau_{-v_i,P}(p)$ is uniformly distributed mod $d$.
\end{lem}
\begin{proof}
    Since the heights of the blocks $P_j$ are $\bar{h}+O(1)$ and $w$ is sufficiently large, the value of $\tau_{-v_i,P}(p)$ is determined by either a single block $P_j$, or two consecutive blocks $P_j,P_{j+1}$ if $\pi_{-v_i,P}(p)$ is close to their common boundary. Write $b=-y_i>0$ and $f(h)=\floor{(h-p_y)/b}$ for the exit time at height $h$.

    In the first case, $\tau_{-v_i,P}(p)=f(h_j)$, where $h_j$ has distribution $\cH$. Since this distribution can be represented as $Z-bU$, where $U$ is uniform on $\{0,\ldots,d-1\}$ and independent of $Z$, $f(Z-bU)$ is uniformly distributed mod $d$.

    In the second case, $\tau_{-v_i,P}(p)$ is determined by $h_j,h_{j+1}$. If $-v_i$ is vertical, then only one block is involved. Otherwise, by symmetry, assume that $-v_i$ points to the right.

    Let $T$ be the largest integer such that $p-Tv_i$ lies to the left of the boundary between $P_j$ and $P_{j+1}$, and let $q$ be its $y$-coordinate. Since $|h_j-h_{j+1}|\leq 1$ and $b\geq 1$, we have
    \[\tau_{-v_i,P}(p)=\begin{cases}
            f(h_j)     & \text{if } h_j<q,             \\
            f(h_j)     & \text{if } h_{j+1}<q\leq h_j, \\
            f(h_{j+1}) & \text{if } q\leq h_j,h_{j+1}.
        \end{cases}\]
    Since $(h_j,h_{j+1})$ is reversible and the chain is in its stationary distribution, the joint distribution of $(h_j,h_{j+1})$ is identical to that of $(h_{j+1},h_j)$. Thus, conditioned on the event that $q\leq h_j,h_{j+1}$, the distribution of $h_{j+1}$ is the same as that of $h_j$. Hence, $\tau_{-v_i,P}(p)$ has the same distribution as $f(h_j)$, which is uniformly distributed mod $d$.
\end{proof}

Next, we show that the $\tau$-distances from a set of points are approximately independent if their projections are far apart.
\begin{lem} \label{lem:cap-tau-indep}
    There exists a constant $c=c_{\ref{lem:cap-tau-indep}}>0$ depending on $w$ but not $\bar{h}$ such that the following holds. Let $(p,i), (p_1,i_1),\ldots,(p_m,i_m)$ be pairs with $i,i_j\in\{2,\ldots,k\}$, all of whose projections $\pi_{-v_i,P}(p),\pi_{-v_{i_j},P}(p_j)$ lie in $[O(1),n-O(1)]$. Suppose there is $\Delta>0$ such that for each $j=1,\ldots,m$,
    \[
        |\pi_{-v_i,P}(p)-\pi_{-v_{i_j},P}(p_j)|\geq \Delta.
    \]
    Let $E$ be any event of positive probability depending on $\tau_{-v_{i_j},P}(p_j)$ for $j=1,\ldots,m$. Then for any $r\in \ZZ$, we have
    \[
        \PP(\tau_{-v_i,P}(p)\equiv r\pmod{d}\mid E)=\frac{1}{d} \pm O(e^{-c\Delta}).
    \]
\end{lem}
\begin{proof}
    Let $J$ be the one or two consecutive block indices determining $\tau_{-v_i,P}(p)$, and define $J_j$ similarly for $\tau_{-v_{i_j},P}(p_j)$. Since the blocks indexed by $J$ are close to the projection $\pi_{-v_i,P}(p)$ (and similarly for the $J_j$), $J$ is at distance at least $\Delta/w-O(1)$ from every other $J_j$. By \Cref{lem:markov-mixing},
    \begin{align*}
        \PP(\tau_{-v_i,P}(p)\equiv r\pmod{d}\mid E)
         & =\PP(\tau_{-v_i,P}(p)\equiv r\pmod{d})\pm O(e^{-c_{\ref{lem:markov-mixing}}(\Delta/w-O(1))}) \\
         & =\frac{1}{d} \pm O(e^{-c\Delta}),
    \end{align*}
    with $c=c_{\ref{lem:markov-mixing}}/w$.
\end{proof}

\subsection{Random properties of the jagged polygon}

\begin{figure}
    \centering
    \begin{tikzpicture}[
            x=2.4cm, y=2.4cm, line cap=round, line join=round,
            polygon boundary/.style={line width=0.75pt, densely dotted},
            jagged boundary/.style={line width=0.95pt},
            ray/.style={->, >=stealth, line width=0.85pt},
            distance measure/.style={<->, >=latex, draw=blue!55!black,
                    line width=0.55pt, dash pattern=on 3pt off 2pt},
            distance guide/.style={draw=blue!40!black, line width=0.35pt},
            every node/.style={font=\large}
        ]
        \def\capHeight{0.245}
        \coordinate (ray start) at (1.15,0.90);
        \coordinate (ray projection) at ({2.30+\capHeight},0.90);
        \coordinate (ray exit) at (2.84,0.90);
        \draw[polygon boundary]
        (0.10,2.50) -- (1.45,2.50) -- (2.30,1.65)
        -- (2.30,0.10) -- (0.95,0.10) -- (0.10,0.95) -- cycle;
        \draw[polygon boundary]
        ({0.10-\capHeight},{2.50+\capHeight})
        -- ({1.45+(sqrt(2)-1)*\capHeight},{2.50+\capHeight})
        -- ({2.30+\capHeight},{1.65+(sqrt(2)-1)*\capHeight})
        -- ({2.30+\capHeight},{0.10-\capHeight})
        -- ({0.95-(sqrt(2)-1)*\capHeight},{0.10-\capHeight})
        -- ({0.10-\capHeight},{0.95-(sqrt(2)-1)*\capHeight})
        -- cycle;

        \draw[jagged boundary]
        (0.10,2.50)
        -- (0.241,2.745) -- (0.338,2.815) -- (0.435,2.705)
        -- (0.532,2.775) -- (0.629,2.825) -- (0.726,2.715)
        -- (0.824,2.795) -- (0.921,2.835) -- (1.018,2.705)
        -- (1.115,2.790) -- (1.212,2.815) -- (1.309,2.745)
        -- (1.45,2.50)
        -- (1.723,2.573) -- (1.754,2.486) -- (1.898,2.512)
        -- (1.883,2.378) -- (2.006,2.383) -- (1.983,2.242)
        -- (2.106,2.247) -- (2.197,2.220) -- (2.175,2.079)
        -- (2.294,2.080) -- (2.286,1.954) -- (2.373,1.923)
        -- (2.30,1.65)
        -- (2.545,1.509) -- (2.620,1.393) -- (2.520,1.278)
        -- (2.585,1.163) -- (2.505,1.048) -- (2.600,1.025)
        -- (2.735,0.975) -- (ray exit) -- (2.700,0.845)
        -- (2.495,0.817) -- (2.575,0.702) -- (2.625,0.587)
        -- (2.510,0.472) -- (2.595,0.357) -- (2.545,0.241)
        -- (2.30,0.10)
        -- (2.159,-0.145) -- (2.062,-0.105) -- (1.965,-0.195)
        -- (1.868,-0.235) -- (1.771,-0.115) -- (1.674,-0.180)
        -- (1.576,-0.100) -- (1.479,-0.220) -- (1.382,-0.125)
        -- (1.285,-0.205) -- (1.188,-0.105) -- (1.091,-0.145)
        -- (0.95,0.10)
        -- (0.677,0.027) -- (0.568,0.036) -- (0.583,0.170)
        -- (0.471,0.176) -- (0.384,0.207) -- (0.410,0.350)
        -- (0.283,0.342) -- (0.281,0.458) -- (0.140,0.436)
        -- (0.124,0.537) -- (0.118,0.649) -- (0.027,0.677)
        -- (0.10,0.95)
        -- (-0.145,1.091) -- (-0.110,1.207) -- (-0.215,1.322)
        -- (-0.180,1.437) -- (-0.095,1.552) -- (-0.200,1.667)
        -- (-0.125,1.783) -- (-0.225,1.898) -- (-0.105,2.013)
        -- (-0.190,2.128) -- (-0.220,2.243) -- (-0.145,2.359)
        -- cycle;

        \node[anchor=south west] at (-0.18,2.86) {$\ol{B}_1^+$};
        \node[anchor=south] at (1.03,2.90) {$B_2$};
        \node[anchor=north] at (0.78,2.43) {$\ol{B}_1$};

        \draw[distance guide] (ray start) -- ($(ray start)+(0,0.36)$);
        \draw[distance guide] (ray exit) -- ($(ray exit)+(0,0.36)$);
        \draw[distance measure]
        ($(ray start)+(0,0.30)$) --
        node[above=3pt, text=blue!55!black, font=\normalsize]
        {$\sim \overline{\tau}_v(p)$}
        ($(ray exit)+(0,0.30)$);

        \draw[ray] (ray start) -- (3.15,0.90);
        \draw[ray] ($(ray start)+(0.16,0)$) --
        node[below=3pt] {$v$} ($(ray start)+(0.55,0)$);
        \fill (ray start) circle[radius=1.6pt] node[below left=3pt] {$p$};
        \fill (ray projection) circle[radius=1.6pt];
        \node[anchor=north west, inner sep=1pt]
        at ($(ray projection)+(0.14,-0.14)$) {$\pi_v(p)$};
    \end{tikzpicture}
    \caption{Illustration of $B_2$, with the inner compact polygon $\ol{B}_1$ and outer compact polygon $\ol{B}_1^+$. The ray from $p$ in direction $v$ meets the outer polygon $\ol{B}_1^+$ at $\pi_v(p)$, while $\overline{\tau}_v(p)$ measures its distance (in steps of $v$) to $\partial B_2$.}
    \label{fig:B2}
\end{figure}
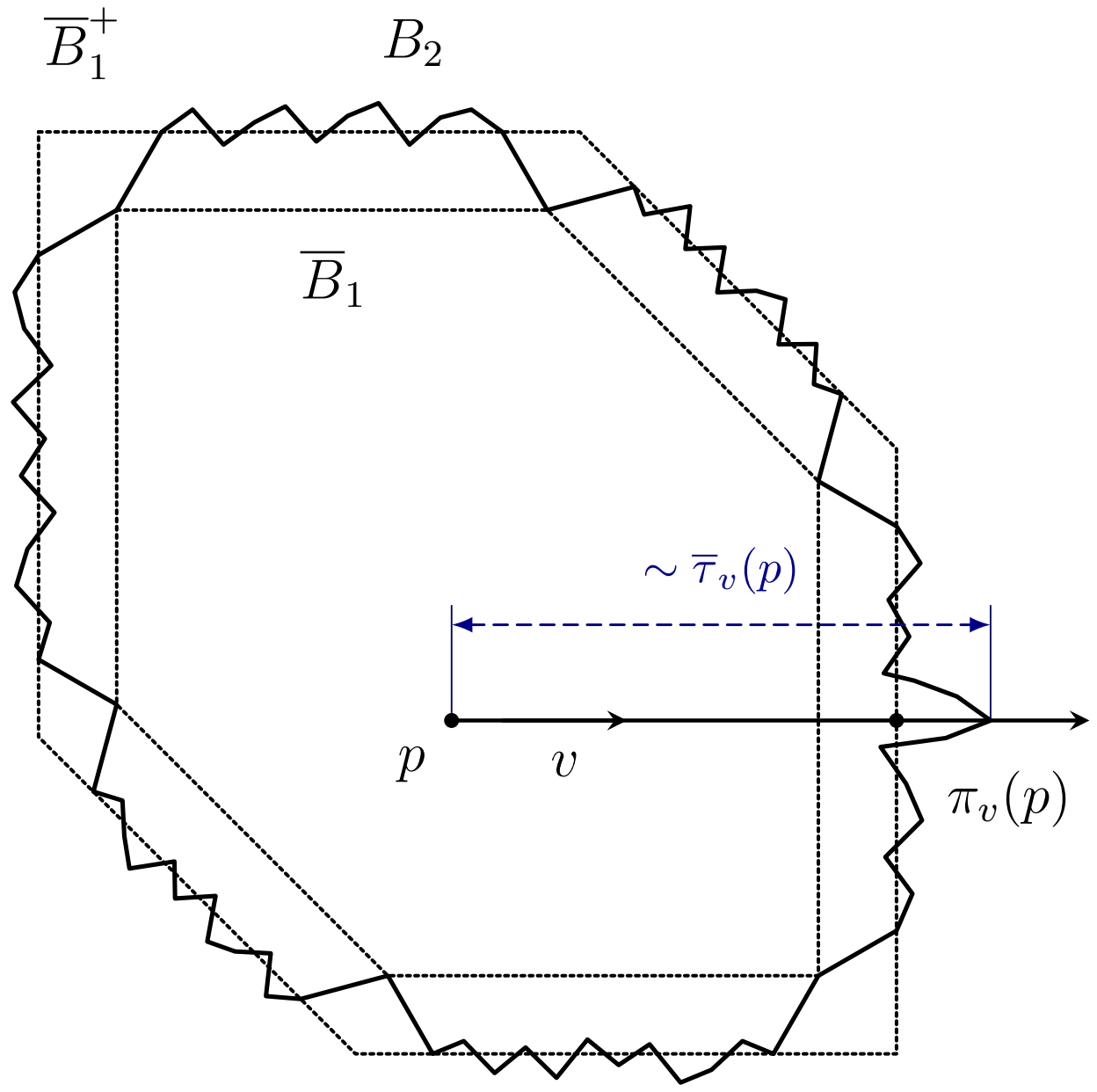

On each side of $B_1$, we add the jagged cap $P$ constructed above, independently. Let $B_2$ be the resulting set. Define $\ol{B}_1^+$ to be the polygon obtained by moving the boundary of $\ol{B}_1$ outward, by the height $\bar{h}$ in the corresponding local coordinates. Thus each side is shifted outwards by distance $\Theta(\bar{h})$, and the boundary of $B_2$ (not too close to a vertex) oscillates only $O(1)$, independently of $\bar{h}$, from the corresponding side of $\ol{B}_1^+$. See \Cref{fig:B2}.

\begin{defn}
    Define a \emph{ray} to be a pair $(p,v)$ where $p\in B_1$ and $v\in \{\pm v_1,\ldots,\pm v_k\}$ is a direction. Let $(p,v)$ be a ray.
    \begin{itemize}
        \item Define the \emph{distance from $p$ to $\partial B_2$ along $v$}, denoted $\overline{\tau}_{v}(p)$, to be the largest integer $t$ such that $p+tv\in B_2$. Denote by $\tau_v(p)$ the image of $\overline{\tau}_{v}(p)$ in $\ZZ/d\ZZ$.

        \item Define the \emph{projection from $p$ to $\partial\ol{B}_1^+$ along $v$}, denoted $\pi_{v}(p)$, to be the intersection of the ray $p+\RR_{\geq 0} v$ with the boundary of $\ol{B}_1^+$. In particular, $\pi_{v}(p)\in \partial \ol{B}_1^+\subset \RR^2$.
    \end{itemize}
\end{defn}

See \Cref{fig:B2} for an illustration of $\overline{\tau}_v(p)$ and $\pi_v(p)$.

Note that in \Cref{lem:cap-tau-unif} and \Cref{lem:cap-tau-indep}, the conclusion holds for a ray $(p,v)$ only if the projection $\pi_{v}(p)$ is far enough away from any corner.

\begin{defn} \label{defn:far-corner}
    Let $C=C_{\ref{defn:far-corner}}>0$ be a constant. We say that a ray $(p,v)$ is \emph{far-corner} if $\pi_{v}(p)$ is at distance at least $C$ from every vertex of $\ol{B}_1^+$. Otherwise, the pair $(p,v)$ is called a \emph{near-corner pair}. A line meeting $B_1$ and parallel to $v$ is a \emph{far-corner line} if both orientations along it are far-corner, and a \emph{near-corner line} otherwise.
\end{defn}

Choose $C_{\ref{defn:far-corner}}$ sufficiently large, depending on the end-piece widths $w_L,w_R$, so that every far-corner ray exits through the random part of a single cap. Then \Cref{lem:cap-tau-unif} and \Cref{lem:cap-tau-indep} give the following.

\begin{lem} \label{lem:B2-unif}
    Let $(p,v)$ be a far-corner ray. Then $\tau_{v}(p)$ is uniformly distributed in $\ZZ/d\ZZ$.
\end{lem}
\begin{proof}
    This is immediate from \Cref{lem:cap-tau-unif}.
\end{proof}

\begin{lem} \label{lem:B2-indep}
    There exist constants $c=c_{\ref{lem:B2-indep}}>0$ and $C=C_{\ref{lem:B2-indep}}>0$ not depending on $\bar{h}$ such that the following holds. Let $(p,u),(p_1,u_1),\ldots,(p_m,u_m)$ be far-corner rays. Suppose that there is $\Delta>0$ such that for each $j=1,\ldots,m$,
    \[\norm{\pi_{u}(p)-\pi_{u_j}(p_j)}\geq \Delta.\]
    Then for any $r\in \ZZ/d\ZZ$ and any event $E$ of positive probability depending only on $\tau_{u_j}(p_j)$ for $j=1,\ldots,m$, we have
    \[\PP(\tau_{u}(p)=r \mid E)=\frac{1}{d} \pm Ce^{-c\Delta}.\]
\end{lem}
\begin{proof}
    The event $\tau_{u}(p)=r$ is independent of all the caps except the one met by $(p,u)$. The projections $\pi_{u_j}(p_j)$ on the same side of $\ol{B}_1^+$ as $\pi_{u}(p)$ are separated from $\pi_{u}(p)$ by $\Omega(\Delta)$, with a factor independent of $\bar{h}$. Then \Cref{lem:cap-tau-indep} gives the desired bound.
\end{proof}

The following useful form of \Cref{lem:B2-indep} allows us to assume independence at the cost of a small probability, even after conditioning on boundary residues at other separated projections.
\begin{lem} \label{lem:B2-coupling}
    Let $(p_1,u_1),\ldots,(p_m,u_m)$ and $(p'_1,u'_1),\ldots,(p'_s,u'_s)$ be far-corner rays, where $s$ may be zero. Suppose that there is $\Delta>0$ such that
    \[
        \norm{\pi_{u_i}(p_i)-\pi_{u_j}(p_j)}\geq\Delta
        \quad\text{for }i\neq j,
    \]
    and
    \[
        \norm{\pi_{u_i}(p_i)-\pi_{u'_j}(p'_j)}\geq\Delta
        \quad\text{for all }i\in[m],\ j\in[s].
    \]
    Let $E$ be any event of positive probability determined by $(\tau_{u'_j}(p'_j))_{j=1}^s$, and write $\PP_E=\PP(\,\cdot\mid E)$. Under $\PP_E$, there is a coupling of $(\tau_{u_i}(p_i))_{i=1}^m$ with independent uniform random variables $X_1,\ldots,X_m\in\ZZ/d\ZZ$ such that
    \[
        \PP_E\paren{(X_i)_{i=1}^m\neq(\tau_{u_i}(p_i))_{i=1}^m}
        =O(md\,e^{-c_{\ref{lem:B2-indep}}\Delta}).
    \]
\end{lem}
\begin{proof}
    Work under $\PP_E$, and write $R_j=\tau_{u_j}(p_j)$ and $c=c_{\ref{lem:B2-indep}}$. We will couple the pairs $(R_j,X_j)$ sequentially. For every $j\in[m]$ and every history $R_l=r_l$ for $l<j$, \Cref{lem:B2-indep}, applied conditional on the event $E\cap\{R_l=r_l\text{ for }l<j\}$, gives
    \[
        \PP_E(R_j=r\mid R_l=r_l\text{ for }l<j)
        =\frac{1}{d}+O(e^{-c\Delta})
    \]
    for every $r\in \ZZ/d\ZZ$. Thus the law of $R_j$, conditional on any previous history $R_l=r_l$ for $l<j$, is at total variation distance $O(de^{-c\Delta})$ from the uniform distribution on $\ZZ/d\ZZ$. Thus, conditional on the previous pairs $(R_l,X_l)$ for $l<j$, we can couple $R_j$ with a uniform $X_j$ with disagreement probability $O(de^{-c\Delta})$. Then the $X_j$ are independent. A union bound therefore gives
    \[
        \PP_E\paren{(R_j)_{j=1}^m\neq(X_j)_{j=1}^m}
        \leq\sum_{j=1}^m\PP_E(R_j\neq X_j)
        =O(md\,e^{-c\Delta}),
    \]
    as required.
\end{proof}

For $\bar{h}$ sufficiently large in terms of $c_{\ref{lem:B2-indep}}$ and $C_{\ref{lem:B2-indep}}$, we have the following.
\begin{lem} \label{lem:B2-pos-point}
    There exists a constant $c=c_{\ref{lem:B2-pos-point}}>0$ such that the following holds. Let $p\in B_1$ and $u_1,\ldots,u_s$ be distinct directions such that $(p,u_j)$ is far-corner for each $j$. Let $(p_1,u'_1),\ldots,(p_m,u'_m)$ be far-corner pairs. Suppose that there is $\Delta>0$ such that for each $i,j$,
    \[
        \norm{\pi_{u_i}(p)-\pi_{u'_j}(p_j)}\geq \Delta.
    \]
    Then for any $r_1,\ldots,r_s\in \ZZ/d\ZZ$ and any event $E$ of positive probability depending only on $\tau_{u'_j}(p_j)$ for $j=1,\ldots,m$, we have
    \[
        \PP(\tau_{u_i}(p)=r_i \text{ for all } i\mid E) > c - O(e^{-c_{\ref{lem:B2-indep}}\Delta}).
    \]
\end{lem}
\begin{proof}
    Since $p$ is at distance at least $\Omega(\bar{h})$ from the boundary of $\ol{B}_1^+$, the projections $\pi_{u_i}(p)$ are pairwise at distance at least $C_3=\Omega(\bar{h})$. Write $c_1=c_{\ref{lem:B2-indep}}$. Apply \Cref{lem:B2-coupling} to these $s\leq 2k$ rays, conditional on $E$, with separation parameter $\min\{\Delta,C_3\}$. For some constant $C_2>0$ independent of $\bar h$, this gives
    \begin{align*}
        \PP(\tau_{u_i}(p)=r_i\text{ for all }i\mid E)
         & \geq d^{-s}-C_2e^{-c_1\min\{\Delta,C_3\}} \\
         & \geq d^{-s}-C_2e^{-c_1C_3}-C_2e^{-c_1\Delta}.
    \end{align*}
    Take $\bar h$ large enough that $C_2e^{-c_1C_3}\leq 1/(2d^{2k})$. Since $s\leq 2k$, the lemma follows with $c=1/(4d^{2k})$.
\end{proof}

Let $\cL_i$ be the set of lines parallel to $v_i$ which intersect $B_1$, and let $\cL=\bigcup_{i=1}^k \cL_i$. Among lines in the prescribed directions, $B_2$ is tileable along every line outside $\cL$. Indeed, a line parallel to $v_i$ missing $B_1$ can only intersect $B_2$ at a cap on a side parallel to $v_i$. From the construction of the caps, the caps on sides parallel to $v_i$ are tileable by $D_i$. For $\ell\in\cL$, the intersection $\ell\cap B_2$ is a contiguous segment, by the construction of the caps. Its length is the number of points in $\ell\cap B_2$.

\begin{defn}
    Denote by $\cT_i(B_2)$ the set of lines $\ell\in \cL_i$ such that $\ell\cap B_2$ is a contiguous segment of length divisible by $d$. Denote by $\cT(B_2)$ the union $\bigcup_{i=1}^k \cT_i(B_2)$. We also call such lines \emph{pure}.

    For a set $B\subseteq B_2$, denote by $\cT_i(B)$ the set of lines $\ell\in \cT_i(B_2)$ such that $\ell\cap B=\ell\cap B_2$. Similarly denote by $\cT(B)$ the union $\bigcup_{i=1}^k \cT_i(B)$. For a point $a\in\ZZ^2$, we write $a\in\cT_i(B)$ to mean that $a$ lies on a line in $\cT_i(B)$.
\end{defn}

In particular, $B$ is tileable along each line in $\cT(B)$. For $B_2$, since $\ell\cap B_2$ is a contiguous segment for each $\ell\in \cL$, we have that $B_2$ is not tileable along each line in $\cL\setminus \cT(B_2)$.

In subsequent sections, we will be constructing sets $B$ by poking holes within $B_1$ in $B_2$. In particular, every such set $B$ will satisfy $B\setminus B_1 = B_2\setminus B_1$.

\subsection{Good translates}

Recall that a direction is one of the vectors $\pm v_1,\ldots,\pm v_k$.

\begin{defn}[Requirements]
    A \emph{directional requirement} is a pair $(u,r)$ where $u$ is a direction and $r\in \ZZ/d\ZZ$. A point $p\in B_1$ \emph{satisfies the directional requirement $(u,r)$} if $\tau_u(p)=r$.

    A \emph{requirement} $R$ is a set of directional requirements $(u_1,r_1),\ldots,(u_m,r_m)$ where $u_i$ are distinct directions and $r_i\in \ZZ/d\ZZ$. We say that a point $p\in B_1$ \emph{satisfies the requirement $R$}, or that \emph{$R(p)$ holds}, if $p$ satisfies each directional requirement. We say that $R$ does not depend on a direction $v_j$ if none of the $u_i$ are $\pm v_j$.
\end{defn}

\begin{rmk}
    The condition that $a\in \cT_i(B_2)$ can be enforced as a requirement. Indeed, the condition $a\in \cT_i(B_2)$ is equivalent to $\tau_{v_i}(a)+\tau_{-v_i}(a)\equiv -1\pmod{d}$, so one could include the pairs $(v_i,r)$ and $(-v_i,-r-1)$ for any $r\in \ZZ/d\ZZ$ in the requirement. If a point $a\in B_1$ satisfies this requirement, then $a\in \cT_i(B_2)$.
\end{rmk}

\begin{lem}[Good point translates exist] \label{lem:good-point-translate}
    Let $v$ be a direction. Let $a\in B_1$ and $R$ be a requirement that does not depend on $v$. Suppose for any positive integer $t\leq \log^4 n$, $a+tdv\in B_1$. Then with probability $1-e^{-\Omega(\log^2 n)}$, there exists a positive integer $t\leq \log^4 n$ such that $R(a+tdv)$ holds.
\end{lem}
\begin{proof}
    For any $b\in B_1$, let $P(b)$ be the set of points $\setcond{\pi_u(b)}{\text{directions }u\neq\pm v}$. Let $\Delta=\log^2 n$. Note that at most $O(1)$ values of $t$ can give a near-corner ray $(a+tdv,u)$ for any direction $u\neq \pm v$.

    Avoid these values of $t$, then pick $m=\Omega(\log^4 n/\Delta)$ positive integers $t_1,\ldots,t_m\leq\log^4 n$ greedily for which all the relevant rays are far-corner and the sets $P(a+t_jdv)$ are pairwise $\Delta$-separated. By \Cref{lem:B2-pos-point},
    \[
        \PP(R(a+t_jdv) \text{ holds}\mid R(a+t_ldv) \text{ does not hold for all } l<j) \geq c_{\ref{lem:B2-pos-point}}-O(e^{-c_{\ref{lem:B2-indep}}\Delta}) > c_{\ref{lem:B2-pos-point}}/2.
    \]
    Thus, the probability that $R(a+tdv)$ fails for every positive integer $t\leq\log^4 n$ is at most $(1-c_{\ref{lem:B2-pos-point}}/2)^m=e^{-\Omega(\log^2 n)}$.
\end{proof}

\begin{defn}[Projection-separated]
    Let $K>0$. We say that a subset $A\subset B_1$ is \emph{$(K,v_i)$-projection-separated} if the following holds. For any pair of points $a,a'\in A$, let $\ell,\ell'$ be the lines parallel to $v_i$ containing $a,a'$ respectively. Then either $\ell=\ell'$ or the distance between $\ell$ and $\ell'$ is at least $K$.
\end{defn}

\begin{defn}
    Let $A\subset \ZZ^2$ and $(R_a)_{a\in A}$ be a family of requirements indexed by $A$. We say that $(R_a)_{a\in A}$ is \emph{compatible} if whenever $a'=a+qu$ for some integer $q$ and $(u,r)\in R_a, (u,r')\in R_{a'}$, we have
    \[r'\equiv r-q\pmod{d}.\]
\end{defn}

Observe that, for a compatible family of requirements with compatible directional requirements $(u,r)\in R_a$ and $(u,r')\in R_{a'}$ on the same line parallel to $u$ and any translation vector $t$, $a+t$ satisfies $(u,r)$ if and only if $a'+t$ satisfies $(u,r')$, whenever both translates lie in $B_1$.

\begin{lem}[Good translates exist] \label{lem:good-translate}
    There exists a constant $c=c_{\ref{lem:good-translate}}>0$ such that the following holds. Let $v=\pm v_i$ be a direction. Let $A\subset B_1$ be of size at most $\log\log\log n$, be within distance $cn$ from the origin, and be $(\Omega(\log^2 n),v_j)$-projection-separated for all $j\neq i$. Let $(R_a)_{a\in A}$ be a compatible family of requirements, none of which depends on $v_i$. Then with probability $1-e^{-\Omega(\log^2 n)}$, there exists a positive integer $t\leq \log^5 n$ such that $R_a(a+tdv)$ holds for all $a\in A$.
\end{lem}
\begin{proof}
    First, for each direction $u$ and line parallel to $u$, if there are multiple points $a\in A$ on this line with directional requirements $(u,r)\in R_a$ for some $r$, then retain just one such directional requirement $(u,r)$ for one of the points $a$ on this line. By compatibility, this does not affect the final conclusion. Now, for each direction $u$, the projections associated with the retained requirements are distinct and $\Omega(\log^2 n)$-separated.

    As in the proof of \Cref{lem:good-point-translate}, let $P(S)$ be the set of points
    \[\setcond{\pi_u(a)}{a\in S, \text{direction } u\neq \pm v_i}\]
    for any set $S\subset B_1$. Since $A$ is $(\Omega(\log^2 n),v_j)$-projection-separated and $A$ is close to the origin, for sufficiently small $c_{\ref{lem:good-translate}}$ the set of projections $P(A)$ is $\Omega(\log^2 n)$-separated. Let $\Delta=\Theta(\log^2 n)$, so that $P(A)$ is $\Delta$-separated.

    By our choice of the polygon $\ol{B}_1$, the projections of the origin in all directions are $\Omega(n)$ from every vertex. The same holds for points within distance $c_{\ref{lem:good-translate}}n$ of the origin when this constant is sufficiently small. The same also holds for $\ol{B}_1^+$ since $\Theta(\bar{h})=o(n)$. Thus all the translated points lie in $B_1$ and all the rays involved in $P(A+tdv)$ are far-corner.

    Using the same greedy argument as in \Cref{lem:good-point-translate}, we can find $m=\Omega(\log^5 n/(\Delta|A|^2))=\Omega(\log^{2.5} n)$ positive integers $t_1,\ldots,t_m\leq\log^5 n$ such that the sets $P(A+t_jdv)$ are pairwise $\Delta$-separated. Note that $P(A+t_jdv)$ itself is $\Delta$-separated, since $P(A)$ is. Write $R(A+u)$ for the assertion that $R_a(a+u)$ holds for all $a\in A$.

    Apply \Cref{lem:B2-coupling} to the at most $2k|A|$ retained directional requirements at each translate, conditional on failure at all earlier translates. This gives
    \[
        \PP(R(A+t_jdv)\text{ holds}\mid R(A+t_ldv)\text{ fails for all }l<j)
        \geq d^{-2k|A|}-O(|A|e^{-c_{\ref{lem:B2-indep}}\Delta})
        \geq c_1/2,
    \]
    where $c_1=d^{-2k|A|}=\Omega(\log^{-0.5}n)$ since $|A|\leq\log\log\log n$. Thus the probability of failure for every positive integer $t\leq\log^5 n$ is at most $(1-c_1/2)^m=e^{-\Omega(\log^2 n)}$.
\end{proof}

There are $(d+1)^{2k}=O(1)$ different requirements. Since $|B_1|=O(n^2)$, the number of sets $A\subset B_1$ of size at most $\log\log\log n$, together with assignments of requirements to their points, is $\exp(O(\log n\,\log\log\log n))=\exp(o(\log^2 n))$. Thus a union bound shows that, with probability $1-e^{-\Omega(\log^2 n)}$, \Cref{lem:good-point-translate} and \Cref{lem:good-translate} hold simultaneously for all choices satisfying their hypotheses. In the rest of the paper, we only invoke \Cref{lem:good-translate} with $|A|=O(1)$. The reason for \Cref{lem:good-translate} having the slowly growing bound $|A|\leq \log\log\log n$ is to make the statement cleaner, without having constants depending on later parts of the proof.

\section{Almost-perfect matching} \label{sec:match}

This section implements Step 3 of the proof strategy.
As discussed above, we will reframe our problem as finding a perfect matching in an auxiliary hypergraph,
which sets up the two-step strategy of (i) finding an almost-perfect matching,
then (ii) applying our new absorption method to resolve the unmatched vertices.
We intend to partition each non-tileable line (in one of our fixed directions $v_1,\dots,v_k$) into segments,
each having a prescribed residue modulo $d$ of a position in the segment where we wish to poke a hole,
so that the remaining punctured line will be tileable. For example, for a line with length $1$ modulo $d$, the number of segments should also be $1$ modulo $d$,
and starting from one end, the residues of hole distances from that end should cycle through $1$ to $d$ repeatedly, starting and ending with $1$.
Once the segments are fixed, the choice of holes corresponds to a perfect matching in the $k$-uniform $k$-partite linear hypergraph $\cH$,
where vertices correspond to segments and edges correspond to cells where poking a hole would satisfy all residue requirements of segments containing it.

A natural first attempt at this strategy is to take all segments to have length about $s$,
with the trade-off that $s$ is small enough that there are many segments in each line,
yet also large enough to benefit from concentration of measure estimates on rays from the segment to the random boundary
(away from the boundary, such distances are uniform modulo $d$, so meet the requirement of any other overlapping segment with probability $1/d$).
For any fixed segment, identified with a vertex $v \in V(\cH)$, its degree in $\cH$ is a sum of indicator random variables, one for each cell in the segment $v$,
mostly having about the same probability of being an allowed edge, and being roughly independent unless they are close together. It is intuitively plausible
that concentration of measure would ensure that $\cH$ is nearly regular and so has an almost-perfect matching by the well-known nibble machinery.

However, a more careful analysis reveals a fundamental problem with this naive strategy.
Firstly, we should clarify what we need from an almost-perfect matching:
the number of segments we fail to cover should be $o(n)$, ideally $O(n^{1-\delta})$ for some $\delta>0$.
If we take $\Delta \sim n^{1-\eps}$ then we have $\Theta(n^{1+\dD})$ segments, of which we need to cover all but an $o(n^{-\dD})$ proportion.
We need to decompose general lengths of size $\Theta(n)$ into segments of length about $\Delta$, which uses $\Theta(n^\eps)$ segments,
so the segment lengths must vary by $\Omega(n^{1-2\eps})$ to cover all possible lengths. However, this means that $\cH$ is only nearly regular
with accuracy $O(n^{-\eps})$, so it is hopeless to cover all but an $o(n^{-\eps})$ proportion of its vertices by a nibble.
At first one might dismiss this as merely a technical inconvenience, but it is surprisingly difficult to circumvent.
Our solution is to randomly vary the lengths of the segments and introduce a weight function on the edges,
chosen to counterbalance the variation in the lengths of segments.

Let $\delta=\frac{1}{200k^2}$, $\varepsilon=\delta/10$ and $\Delta=n^{1-\varepsilon}$.

For each line $\ell\in \cL$, say $\ell\in \cL_i$, let $n_\ell=|\ell\cap B_1|$, so $n_\ell=\Theta(n)$. The construction of $B_2$ ensures that $\ell\cap B_2$ is also a contiguous segment, of length $\overline{\tau}_{v_i}(x)+\overline{\tau}_{-v_i}(x)+1$ for any $x\in \ell\cap B_1$. Suppose this length is $r_\ell$ mod $d$, with $r_\ell\in \ZZ/d\ZZ$.

Pick real numbers $E_\ell^+\in [\Delta,2d\Delta]$ and $E_\ell^-\in [-2d\Delta,-\Delta]$ (depending only on $n_\ell$ and $r_\ell$) such that
\[
    \frac{n_\ell+E_\ell^+}{\Delta} \quad \text{ and }\quad \frac{n_\ell+E_\ell^-}{\Delta}
\]
are integers congruent to $r_\ell\pmod{d}$. Let $p_\ell^{\pm}:=1+E_\ell^{\pm}/n_\ell$. Conditional on $B_2$, choose $p_\ell\in \{p_\ell^+, p_\ell^-\}$ independently for different $\ell$, with probabilities such that $\EE[p_\ell\mid B_2]=1$. The conditional distribution of $p_\ell$ depends on the boundary of $B_2$ only through $r_\ell$. Then $p_\ell$ satisfies the following properties.
\begin{enumerate}
    \item $1-O(\Delta/n)< p_\ell^-\leq 1\leq p_\ell^+<1+O(\Delta/n)$.
    \item $\EE[p_\ell\mid r_\ell]=\EE[p_\ell]=1$.
    \item $\frac{p_\ell n_\ell}{\Delta}\equiv r_\ell \pmod{d}$.
\end{enumerate}

Partition $\ell\cap B_1$ into $t_\ell:=\frac{p_\ell n_\ell}{\Delta}$ segments of length $\floor{\Delta/p_\ell}$ or $\ceil{\Delta/p_\ell}$, using a fixed deterministic rule depending only on $n_\ell$ and $p_\ell$. Let these segments be $S_1,\ldots,S_{t_\ell}$ in the order of direction $v_i$. We say that a cell $x$ in $S_j$ is \emph{permissible} if $\tau_{-v_i}(x)\equiv j-1\pmod{d}$. That is, one in every $d$ cells in $S_j$ is permissible. The permissible cells are defined precisely so that if we poke a hole in each $S_j$ at a permissible cell, then the resulting subset of $\ell\cap B_2$ is tileable by $D_i$.

We are now ready to define the $k$-partite weighted hypergraph $(\cH,\omega)$. Its vertex set has $k$ parts, with the $j$-th part consisting of all segments in the partition of $\ell\cap B_1$ for all $\ell\in \cL_j\setminus \cT_j(B_2)$. For each $x\in B_1$, if $x$ belongs to $k$ segments $S_1,\ldots,S_k$ from different parts in $V(\cH)$, and $x$ is permissible in each $S_i$, then we add the edge $e=(S_1,\ldots,S_k)$ with weight $\omega(e)=\prod_{\ell} p_{\ell}$, where the product is over all lines $\ell\in \cL$ containing $x$. The hypergraph is linear because two segments in different directions intersect in at most one cell, and $1/2<\omega(e)\leq 2$ for sufficiently large $n$.

If we were able to find a perfect matching in $\cH$, then we would be done: poking holes in the common cells corresponding to the edges in the perfect matching would make $B_2$ tileable along every line in $\cL$. We are unable to find a perfect matching directly, but nevertheless, we can find an almost-perfect matching, leaving only $O(n^{1-\delta})$ vertices unmatched. This would leave only $O(n^{1-\delta})$ lines not tileable, which will be fixed later.

\begin{lem} \label{lem:approx-regular}
    With probability $1-o(1)$, the weighted hypergraph $(\cH,\omega)$ defined above is approximately $\frac{(d-1)^{k-1}}{d^{2k-1}}\Delta$-regular in the following sense: all but $O(n^{0.1})$ vertices $S\in V(\cH)$ satisfy
    \[\deg_{\omega}(S):=\sum_{X\ni S} \omega(X) = \frac{(d-1)^{k-1}}{d^{2k-1}}\Delta+O(\sqrt{n}).\]
\end{lem}

Let us briefly describe the heuristic for why the degree should be around $\frac{(d-1)^{k-1}}{d^{2k-1}}\Delta$. A segment $S$ in $V(\cH)$, say parallel to the direction $v_1$, has length approximately $\Delta$, and approximately $\Delta/d$ of its cells are permissible in $S$. For such a cell to correspond to an edge in $\cH$, it must also be permissible in the other $k-1$ segments containing it, with each happening with probability $1/d$. Furthermore, the cell must not lie in $\cT_j(B_2)$ for each $j=2,\ldots,k$, which happens with probability $(d-1)/d$ for each $j$. If all these events are independent, then the expected degree of $S$ would be $\frac{\Delta}{d}\paren{\frac{d-1}{d^2}}^{k-1}$.

Before proving the lemma, we record a standard martingale concentration inequality that we will use.
\begin{thm}[Azuma-Hoeffding inequality~{\cite[Theorem~5.2]{CL06}}] \label{thm:azuma}
    Let $X_0,X_1,\ldots,X_m$ be a martingale such that for each $1\leq i \leq m$, $|X_i - X_{i-1}|\leq c_i$ almost surely. Then for any $t>0$,
    \[\PP(|X_m - X_0| \geq t) \leq 2\exp\left(\frac{-t^2}{2\sum_{i=1}^m c_i^2}\right).\]
\end{thm}

\begin{proof}[Proof of \Cref{lem:approx-regular}]
    Let $D_1=\log^2 n$. We say that a line $\ell$ is \emph{bad} if any of the following holds.
    \begin{itemize}
        \item $\ell$ is a near-corner line.
        \item $\ell$ is within distance $D_1$ of a side of $\ol{B}_1$ parallel to $\ell$.
    \end{itemize}
    We say that a segment $S\in V(\cH)$ is \emph{bad} if the line containing it is bad. There are only at most $O(D_1 n/\Delta)=O(n^{0.1})$ bad vertices in total, so it suffices to show approximate regularity for vertices which are not bad.

    Fix a line $\ell\in\cL_1$ that is not bad, and fix a point $a\in\ell\cap B_1$. Condition on $\tau_{v_1}(a)$, $\tau_{-v_1}(a)$ and $p_\ell$, with $r_\ell\neq 0$.

    Let $S$ be any segment in this fixed partition. Its degree is
    \[\deg_{\omega}(S)=p_\ell \sum_{x\in P}\prod_{j=2}^k p_{\ell_j(x)}\one(E_j(x)\cap F_j(x))=p_\ell \sum_{x\in P} q_x,\]
    where the following definitions apply.
    \begin{itemize}
        \item $P$ is the set of permissible cells in $S$.
        \item $q_x$ is the random variable $\prod_{j=2}^k p_{\ell_j(x)}\one(E_j(x)\cap F_j(x))$.
        \item $\ell_j(x)$ is the line parallel to $v_j$ passing through $x$.
        \item $E_j(x)$ is the event that $\ell_j(x)\notin \cT_j(B_2)$.
        \item $F_j(x)$ is the event that $x$ is permissible in the segment of the partition of $\ell_j(x)\cap B_1$ containing $x$.
    \end{itemize}
    We have $|P|=|S|/d+O(1)=\Delta/(p_\ell d)+O(1)$. Let $P'\subseteq P$ be the subset of points which are at distance at least $D_1$ from the boundary of $\ol{B}_1$ and do not lie on any near-corner line. Since $S$ is not bad, we have $|P\setminus P'|=O(D_1)=O(\log^2 n)$. Also $0\leq q_x\leq 2$ for sufficiently large $n$. Thus, the contribution to $\deg_{\omega}(S)$ from $P\setminus P'$ is at most $O(\log^2 n)$, well within our desired error bound of $O(\sqrt{n})$, so it suffices to show concentration for the sum $\sum_{x\in P'} q_x$.

    We will first compute the expectation. For each $x\in P'$, the variable $q_x$ depends only on $\tau_{\pm v_j}(x)$ and $p_{\ell_j(x)}$ for $j=2,\ldots,k$. Fix $x\in P'$, and write $\ell_j,E_j,F_j$ instead of $\ell_j(x),E_j(x),F_j(x)$ for brevity.

    First consider the ideal case in which $\tau_{\pm v_j}(x)$ are independent uniform random variables in $\ZZ/d\ZZ$ for $j=2,\ldots,k$. Conditional on these residues, choose the $p_{\ell_j}$ independently according to their original distributions determined by $r_{\ell_j}$. Denote probability and expectation under this ideal case by $\PP_0$ and $\EE_0$. Since $r_{\ell_j}=\tau_{v_j}(x)+\tau_{-v_j}(x)+1$, it is uniform in $\ZZ/d\ZZ$. The event $E_j$ occurs if and only if $r_{\ell_j}\neq 0$. Thus $\PP_0(E_j)=\frac{d-1}{d}$. Conditional on $(r_{\ell_j},p_{\ell_j})$, each of $\tau_{v_j}(x)$ and $\tau_{-v_j}(x)$ is still uniform in $\ZZ/d\ZZ$, although they are no longer independent. The partition of $\ell_j\cap B_1$ is now fixed, so $F_j$ prescribes a single residue for $\tau_{-v_j}(x)$, giving $\PP_0(F_j\mid r_{\ell_j},p_{\ell_j})=1/d$. Since $\EE_0[p_{\ell_j}\mid r_{\ell_j}]=1$, we have
    \begin{align*}
        \EE_0[p_{\ell_j}\one(E_j\cap F_j)]
         & = \EE_0[p_{\ell_j}\one(r_{\ell_j}\neq 0)\PP_0(F_j\mid r_{\ell_j}, p_{\ell_j})] \\
         & = \frac{1}{d}\EE_0[p_{\ell_j}\one(r_{\ell_j}\neq 0)] = \frac{d-1}{d^2}.
    \end{align*}
    Thus,
    \[\EE_0[q_x] = \prod_{j=2}^k \EE_0[p_{\ell_j}\one(E_j\cap F_j)] = \paren{\frac{d-1}{d^2}}^{k-1}=:\mu.\]

    In reality, the variables $\tau_{\pm v_j}(x)$ are not independent. Choose $D_2=c_0D_1$ for a sufficiently small fixed constant $c_0>0$. Since $x$ is at distance at least $D_1$ from the boundary, the projections $\pi_{\pm v_j}(x)$ for $j=1,\ldots,k$ are pairwise at distance at least $D_2$. By \Cref{lem:B2-coupling}, conditional on the two fixed residues on $\ell$, we can couple $(\tau_{\pm v_j}(x))_{j=2}^k$ with independent uniform residues with disagreement probability $O(e^{-cD_2})$, where $c=c_{\ref{lem:B2-indep}}$. Comparing the actual and ideal versions (above) of $q_x$ and using $0\leq q_x\leq 2$, we have
    \[\EE[q_x] = \mu + O(e^{-cD_2}).\]
    Since $|P'|=\frac{\Delta}{p_\ell d}+O(\log^2 n)$, we have
    \[\EE[\deg_{\omega}(S)] = \frac{(d-1)^{k-1}}{d^{2k-1}}\Delta+O(\Delta e^{-cD_2}+\log^2 n).\]

    To show concentration, we will partition $P'$ into a small number of subsets $P_1,\ldots,P_m$ such that the summands within each $P_i$ are approximately independent. Then we apply the martingale concentration inequality to the sums over each $P_i$.

    We construct a graph $G$ with vertex set $P'$ and edges as follows. For each $x\in P'$, denote by $Q_x$ the set $\setcond{\pi_{\pm v_j}(x)}{j=2,\ldots,k}$. For $x,x'\in P'$, we add an edge $xx'$ in $G$ if $\norm{q-q'}\leq D_2$ for some $q\in Q_x$ and $q'\in Q_{x'}$. Observe that for each pair of directions $\pm v_j, \pm v_{j'}$ and for each $x$, there are only $O(D_2)$ many points $x'\in \ell\cap B_1$ such that
    \[
        \norm{\pi_{\pm v_j}(x)-\pi_{\pm v_{j'}}(x')}\leq D_2.
    \]
    Thus $G$ has maximum degree $O(D_2)$, so a simple greedy colouring shows that the chromatic number of $G$ is at most $O(D_2)$. Hence, we can partition $P'$ into $m=O(D_2)$ subsets $P_1,\ldots,P_m$ such that each $P_i$ is an independent set in $G$.

    We now show concentration for the sum $\sum_{x\in P_h} q_x$ for each class $P_h$. Focus on the class $P_1$; the argument for the other classes is identical. Let $P_1=\{x_1,\ldots,x_s\}$ and let $\cF_i=\sigma(q_{x_1},\ldots,q_{x_i})$ be the $\sigma$-algebra generated by the random variables $q_{x_1},\ldots,q_{x_i}$. First condition on the boundary residues $\tau_{\pm v_j}(x_t)$ and the values of $p_{\ell_j(x_t)}$ for $t<i$ and $j=2,\ldots,k$. Then the values of $q_{x_1},\ldots,q_{x_{i-1}}$ are determined. Since $P_1$ is an independent set, the projections in the sets $Q_{x_i}$ are pairwise at distance at least $D_2$. Thus \Cref{lem:B2-coupling} couples the new residues, conditional on the exposed boundary residues, with independent uniform ones with error $O(e^{-cD_2})$. The previous calculation then gives
    \[
        \EE[q_{x_i}\mid\cF_{i-1}]=\mu+O(e^{-cD_2}),
        \qquad
        \abs{\EE[q_{x_i}\mid\cF_{i-1}]-\EE[q_{x_i}]}=O(e^{-cD_2}).
    \]
    Define
    \[
        Z_i=q_{x_i}-\EE[q_{x_i}\mid\cF_{i-1}],
        \qquad Y_i=\sum_{j=1}^i Z_j,\qquad Y_0=0.
    \]
    Then $(Y_i)_{i=0}^s$ is a martingale with $|Z_i|\leq 2$, so \Cref{thm:azuma} gives
    \[\PP\paren{|Y_s|\geq t}\leq 2\exp\paren{-\frac{t^2}{8s}}.\]
    Since $q_{x_i}-\EE[q_{x_i}]=Z_i+O(e^{-cD_2})$, there is a constant $C>0$ such that
    \[
        \PP\paren{\abs{\sum_{x\in P_1}q_x-\EE\sqb{\sum_{x\in P_1}q_x}}
            \geq t+C|P_1|e^{-cD_2}}
        \leq 2\exp\paren{-\frac{t^2}{8|P_1|}}.
    \]
    Apply this bound to each class $P_h$ with threshold $t/m$. The triangle inequality and a union bound give
    \begin{align*}
         & \PP\paren{\abs{\sum_{x\in P'}q_x-\EE\sqb{\sum_{x\in P'}q_x}}
                  \geq t+C|P'|e^{-cD_2}} \\
         & \hspace{3em}\leq 2m\exp\paren{-\frac{t^2}{8m^2|P'|}}.
    \end{align*}
    Set $t=\sqrt n$. Since $m=O(\log^2 n)$ and $|P'|=O(n^{1-\varepsilon})$, the last bound is
    \[
        O(\log^2 n)\exp\paren{-\Omega\paren{\frac{n^\varepsilon}{\log^4 n}}}
        =O\paren{e^{-\Omega(\log^2 n)}}.
    \]
    Also $|P'|e^{-cD_2}=o(\sqrt n)$. Multiplying by the fixed bounded factor $p_\ell$ and adding back the contribution of $P\setminus P'$ therefore gives
    \[
        \PP(|\deg_{\omega}(S)-\EE[\deg_{\omega}(S)]|\geq C'\sqrt n)
        =O\paren{e^{-\Omega(\log^2 n)}}
    \]
    for some constant $C'>0$. Taking a union bound over the $O(n^{1+\varepsilon})$ line and segment indices, we obtain with probability $1-o(1)$ that every segment $S$ which is not bad satisfies
    \[
        \deg_{\omega}(S)= \frac{(d-1)^{k-1}}{d^{2k-1}}\Delta+O(\sqrt{n}).
    \]

\end{proof}

We will use the following quantitative result due to Ehard, Glock and Joos~\cite{EGJ20} on almost-perfect matchings in almost-regular hypergraphs.
\begin{thm}[{\cite[Theorem~1.2]{EGJ20}}] \label{thm:almost-matching}
    Let $\delta\in (0,1)$ and $r\in \NN$ with $r\geq 2$, and let $\varepsilon:=\frac{\delta}{50r^2}$. Then there exists $\Delta_0$ such that for all $\Delta\geq \Delta_0$, the following holds. Let $\cH$ be an $r$-uniform hypergraph with $\Delta(\cH)\leq \Delta$ and maximum codegree at most $\Delta^{1-\delta}$ as well as $e(\cH)\leq \exp(\Delta^{\varepsilon^2})$. Suppose that $\cW$ is a set of at most $\exp(\Delta^{\varepsilon^2})$ non-negative weight functions on $E(\cH)$. Then, there exists a matching $\cM$ in $\cH$ such that $\omega(\cM)=(1\pm \Delta^{-\varepsilon})\omega(E(\cH))/\Delta$ for all $\omega\in \cW$ with $\omega(E(\cH))\geq \max_{e\in E(\cH)} \omega(e)\Delta^{1+\delta}$.
\end{thm}

In particular, we will be using the following special case.
\begin{lem} \label{lem:almost-matching-linear}
    Let $\eta=\frac{1}{100k^2}$ and $\Delta>N^{2/3}$. Let $(\cH,\omega)$ be a weighted $k$-uniform $k$-partite linear hypergraph with $N=|V(\cH)|$ and weights satisfying $0\leq \omega(e)\leq 2$ for all $e\in E(\cH)$. Suppose that $(\cH,\omega)$ is approximately $\Delta$-regular in the following sense. For all but $O(N^{1/2})$ vertices $v\in V(\cH)$, we have
    \begin{equation} \label{eq:approx-reg}
        \deg_\omega(v)=\sum_{e\ni v} \omega(e) = \Delta + O(N^{3/5}).
    \end{equation}
    Then there exists a matching $\cM$ in $\cH$ covering all but at most $O(N^{1-2\eta/3})$ vertices of $\cH$.
\end{lem}
\begin{proof}
    Let $V^*$ be the set of vertices satisfying \eqref{eq:approx-reg}, so that $|V(\cH)\setminus V^*|=O(N^{1/2})$. By linearity of $\cH$, we have $\Delta=O(N)$. Let $\cH'$ be the (unweighted) subhypergraph of $\cH$ obtained by keeping each edge $e\in E(\cH)$ with probability $\omega(e)/2$ independently. Then
    \[
        \EE[\deg_{\cH'}(v)] = \frac12 \deg_{\omega}(v).
    \]
    By Hoeffding's inequality, we have for any $t>0$,
    \[
        \PP\paren{|\deg_{\cH'}(v) - \EE[\deg_{\cH'}(v)]| \geq t} \leq 2\exp\paren{\frac{-2t^2}{\deg_{\cH}(v)}}.
    \]
    Setting $t=N^{3/5}$ and using $\deg_{\cH}(v)=O(N)$ gives that the above failure probability is at most $2e^{-\Omega(N^{1/5})}$ for each $v\in V^*$. A union bound therefore shows that, with probability $1-o(1)$, every $v\in V^*$ satisfies
    \[
        \deg_{\cH'}(v) = \frac12 \Delta + O(N^{3/5}).
    \]
    Fix such a choice of $\cH'$. Let $\cH''$ be the restriction of $\cH'$ to the vertices in $V^*$. Since $\cH$ is linear, we have $\deg_{\cH''}(v)\geq \deg_{\cH'}(v)-|V(\cH)\setminus V^*|$. Thus for every $v\in V^*$, we have
    \[
        \deg_{\cH''}(v) = \frac12 \Delta + O(N^{3/5}).
    \]
    Let $\Delta'$ be the maximum degree of $\cH''$. Then $\Delta'=\frac12\Delta+O(N^{3/5})=\Theta(\Delta)$, since $\Delta>N^{2/3}$. Counting edges by summing degrees, we obtain
    \[
        e(\cH'')=\frac{1}{k}(N-O(N^{1/2}))(\Delta'-O(N^{3/5}))
        =\frac{N}{k}(\Delta'-O(N^{3/5})).
    \]
    Apply \Cref{thm:almost-matching} to $\cH''$, with parameter $\delta=1/2$ and the single weight function identically equal to $1$. Then there is a matching $\cM$ in $\cH''$ (and hence in $\cH$) of size at least
    \begin{align*}
        (1-(\Delta')^{-\eta})\frac{e(\cH'')}{\Delta'}
         & =\frac{N}{k}\bigl(1-O(\Delta^{-\eta}+N^{3/5}\Delta^{-1})\bigr) \\
         & =\frac{N}{k}-O(N^{1-2\eta/3}).
    \end{align*}
    Thus $\cM$ covers all but at most $O(N^{1-2\eta/3})$ vertices of $\cH$.
\end{proof}

Another simple application of \Cref{lem:B2-coupling} shows that approximately $1/d$ of the lines in each direction are pure (belonging to $\cT(B_2)$), and approximately $(d-1)/d$ are non-pure. Thus our weighted hypergraph $(\cH,\omega)$ has $N=|V(\cH)|=\Theta(n^{1+\varepsilon})$ with probability $1-o(1)$.

By \Cref{lem:approx-regular} and \Cref{lem:almost-matching-linear}, with $\eta=2\delta$, we obtain a matching $\cM$ leaving at most
\[
    O(N^{1-4\delta/3})
    =O\bigl(n^{(1+\varepsilon)(1-4\delta/3)}\bigr)
    =O(n^{1-\delta})
\]
vertices uncovered, since $\varepsilon=\delta/10$. The matching corresponds to a set of cells $A_3\subset B_1$. Let $\cS$ be the set of segments in $V(\cH)$ not covered by $\cM$. Set $B_3=B_2\setminus A_3$, and we record the main properties of $B_3$ and $\cS$ below.
\begin{enumerate}
    \item $\cS$ consists of $O(n^{1-\delta})$ segments, each of length $\Omega(n^{1-\varepsilon})$.
    \item $B_3$ is tileable along every line $\ell\in \cL$ not containing any segment in $\cS$.
    \item For every line $\ell\in \cL$ containing some segment in $\cS$, $B_3$ can be made tileable along $\ell$ by poking a hole in any permissible cell in each segment in $\cS$ contained in $\ell$.
    \item $\cT(B_2)=\cT(B_3)$.
\end{enumerate}

For the last property, every line in $\cT(B_2)$ is untouched by $A_3$, so it remains in $\cT(B_3)$.

\section{Cover down} \label{sec:cover-down}

Now we will implement Step 4, in which we prepare for absorption by ensuring a more convenient form for the set of untileable lines.
So far we have constructed $B_3$, which is tileable along almost every line in $\cL$. We will poke a small number of holes in $B_3$ to obtain $B_4$, which is tileable along every line in $\cL$ except for a small number of lines which are pure and within $\alpha n$ of the center of $B_3$ for some small $\alpha>0$. There are two main reasons for doing this step. First, by concentrating the untileable lines near the center, we ensure that the holes we poke in the later steps are also near the center (in particular, they are within $B_1$). The second reason is that for any two points $x,x'$ near the center and two distinct directions $v,v'$, the projections $\pi_v(x)$ and $\pi_{v'}(x')$ on the boundary are far from each other. This ensures that the distances $\tau_v(x)$ and $\tau_{v'}(x')$ are essentially independent, allowing us to apply \Cref{lem:good-translate} in the later steps.

Let $\delta=\frac{1}{200k^2}$ and $\varepsilon=\delta/10$ as in the previous section.
\begin{defn}[State]
    A \emph{state} is a pair $(B,\cS)$, where $B\subseteq B_3$ and $\cS$ is a set of segments in $B\cap B_1$, each parallel to a direction and with a prescribed residue class modulo $d$ of permissible cells, satisfying the following.
    \begin{enumerate}
        \item $|B_3\setminus B|=O(n^{1-3\varepsilon})$.
        \item Each segment in $\cS$ has length $\Theta(n^{1-\varepsilon})$. Segments on the same line are disjoint, and no line in $\cT(B)$ contains a segment in $\cS$.
        \item $B$ is tileable along every line $\ell\in \cL$ not containing any segment in $\cS$.
        \item For every line $\ell\in \cL$ containing such a segment, $B$ can be made tileable along $\ell$ by poking a hole in any permissible cell in each segment in $\cS$ contained in $\ell$.
    \end{enumerate}
\end{defn}

Recall that from the previous section, we have a state $(B_3,\cS)$ of $O(n^{1-\delta})$ segments. Let $O\in B_1$ be the origin so that the ball with center $O$ and radius $\alpha n$ is contained in $\ol{B}_1$, where $\alpha >0$ is a small constant. Our aim is to poke holes in $B_3$ and modify $\cS$, eventually reaching a state $(B_4,\cS')$ such that every segment in $\cS'$ is inside the ball with center $O$ and radius $\alpha n$.

\begin{lem} \label{lem:geom-mono}
    There exists a constant $0<c_{\ref{lem:geom-mono}}<1$ such that the following holds. Let $v,v'$ be two non-parallel directions and $X,A\in \RR^2$ be two points. Let $\ell,\ell'$ be the lines through $A$ parallel to $v,v'$ respectively, and let $P,P'$ be the orthogonal projections of $X$ onto $\ell,\ell'$ respectively. Then
    \[|XP'| \leq c_{\ref{lem:geom-mono}} (|XP|+2|AP|),\]
    where $|XP|$ denotes the length of the segment $XP$.
\end{lem}

\begin{figure}
    \centering
    \begin{tikzpicture}[scale=1]
        \coordinate (A) at (0, 0);
        \coordinate (X) at (1, 5);

        \coordinate (P) at (1, 0);

        \coordinate (Pprime) at (3, 3);

        \draw[line width=0.5pt, blue] (-2, 0) -- (4, 0) node[right] {$\ell$};
        \draw[line width=0.5pt, red] (-1, -1) -- (4, 4) node[above right] {$\ell'$};
        \draw[line width=0.5pt, dashed] (X) -- (P);
        \draw[line width=0.5pt, dashed] (X) -- (Pprime);

        \fill[black] (X) circle (2pt) node[above left] {$X$};
        \fill[black] (A) circle (2pt) node[below right] {$A$};
        \fill[blue] (P) circle (2pt) node[below right] {$P$};
        \fill[red] (Pprime) circle (2pt) node[below right] {$P'$};
    \end{tikzpicture}
    \caption{Illustration of \Cref{lem:geom-mono}.}
\end{figure}

\begin{proof}
    Let $\theta^*\in (0,\pi/2]$ be the smallest angle between any two non-parallel directions and set $c=c_{\ref{lem:geom-mono}}=\max(1/2, \cos \theta^*)<1$. Let $\theta$ be the angle between $v$ and $v'$. After a change of coordinates, we can assume that $A$ is the origin, $\ell$ is the $x$-axis, and $\ell'$ has unit direction vector $(\cos \theta, \sin \theta)$. Let $X=(x,y)$. Then $P=(x,0)$ and the length of $XP'$ is $|y\cos\theta-x\sin\theta|$. Thus,
    \begin{align*}
        |XP'| & = |y\cos\theta-x\sin\theta| \leq |y\cos\theta| + |x\sin\theta| \\
              & \leq c(|y|+2|x|)=c(|XP|+2|AP|).
    \end{align*}
\end{proof}

\begin{defn}
    Given a segment $S$, we define its \emph{cost}, denoted by $\cost(S)$, as follows. Let $\ell$ be the line containing $S$ and let $M$ be the midpoint of $S$. Let $P$ be the projection from $O$ onto $\ell$. Then we define
    \[\cost(S) = |OP|+2|MP|.\]
\end{defn}
Note that if $\cost(S)<\alpha n$, then $S$ is contained in the ball with center $O$ and radius $\alpha n+O(n^{1-\varepsilon})$.

\begin{lem} \label{lem:coverdown-step}
    Let $(B,\cS)$ be a state, and $S\in \cS$ with $\cost(S)=s$. Then there exists a state $(B',\cS')$ such that $B'$ is obtained from $B$ by poking a hole in $S$, and $\cS'$ is obtained from $\cS$ by replacing $S$ with at most $dk$ segments, each of cost at most $c_{\ref{lem:geom-mono}}s+O(n^{1-\varepsilon})$.
\end{lem}
\begin{proof}
    Suppose without loss of generality that $S$ lies on the line $\ell$ parallel to $v_1$. Let $M$ be the midpoint of $S$, and choose the sign of each $v_i$, for $i=2,\ldots,k$, so that $M$ moves closer to $O$ in the direction of each $v_i$, i.e. $\ang{M,v_i}\leq 0$. Since the length of $S$ is $\Omega(n^{1-\varepsilon})$, by \Cref{lem:good-point-translate}, there exist at least $\Omega(n^{1-2\varepsilon})$ choices of permissible cells $Y$ in $S$ satisfying the following.
    \begin{enumerate}
        \item For each $i=2,\ldots,k$, the line $\ell_i$ parallel to $v_i$ passing through $Y$ is in $\cT_i(B_3)$.
        \item For each $i=2,\ldots,k$, $\tau_{-v_i}(Y)=0$.
    \end{enumerate}
    Since $B$ differs from $B_3$ by only $O(n^{1-3\varepsilon})$ holes, $\cT(B_3)$ differs from $\cT(B)$ by only $O(n^{1-3\varepsilon})$ lines, so
    there are still $\Omega(n^{1-2\varepsilon})$ choices of $Y$ satisfying the above two conditions for $B$ in place of $B_3$. Pick one such $Y$.

    For $i=2,\ldots,k$, let $P_i$ be the projection of $O$ onto $\ell_i$. Note that $\ell_i\cap B_1$ is a segment of length $\Omega(n)$ on each side of $P_i$. Introduce $d-1$ pairwise disjoint new segments $S_i^1,\ldots,S_i^{d-1}$ in $B_1\cap\ell_i$, each of length $\floor{n^{1-\varepsilon}}$, ordered in direction $v_i$, lying strictly in that direction from $Y$, and with midpoints at distance $O(n^{1-\varepsilon})$ from $P_i$. Since $\ell_i\in\cT_i(B)$, these segments lie in $B$ and $\ell_i$ contains no segment in $\cS$. Declare $Z\in S_i^j$ permissible when $\tau_{-v_i}(Z)\equiv j\pmod d$. Since $\tau_{-v_i}(Y)=0$ and $\ell_i$ is pure, poking holes at $Y$ and at any permissible cell in each $S_i^j$ leaves only intervals of length divisible by $d$ along $\ell_i$. Take $B' = B\setminus \{Y\}$ and $\cS'$ to be the set obtained from $\cS$ by replacing $S$ with $S_i^j$ for $i=2,\ldots,k$ and $j=1,\ldots,d-1$. Then $(B',\cS')$ is a state, so it suffices to show that $\cost(S_i^j) \leq c_{\ref{lem:geom-mono}}\cost(S)+O(n^{1-\varepsilon})$ for each $i,j$.

    If $\ell_i'$ is the line through $M$ parallel to $v_i$ and $P,P'$ are the projections from $O$ onto $\ell, \ell_i'$, then by \Cref{lem:geom-mono}, we have $|OP'| \leq c_{\ref{lem:geom-mono}}(|OP|+2|MP|)$. Since $Y$ is at distance at most $O(n^{1-\varepsilon})$ from $M$ and $P'$ is at distance at most $O(n^{1-\varepsilon})$ from the midpoint of $S_i^j$, we have $\cost(S_i^j)=|OP'|+O(n^{1-\varepsilon})$ and $\cost(S)=|OP|+2|MP|+O(n^{1-\varepsilon})$. Thus, $\cost(S_i^j) \leq c_{\ref{lem:geom-mono}}(|OP|+2|MP|)+O(n^{1-\varepsilon}) = c_{\ref{lem:geom-mono}}\cost(S)+O(n^{1-\varepsilon})$.
\end{proof}

Repeatedly apply \Cref{lem:coverdown-step} whenever there is a segment with cost at least $\alpha n/2$. Put $q=(1+c_{\ref{lem:geom-mono}})/2<1$. At each step, a segment of cost $s\geq \alpha n/2$ is replaced by at most $dk$ segments, each of cost at most $c_{\ref{lem:geom-mono}}s+O(n^{1-\varepsilon})\leq qs$. Since the initial cost of each segment is $O(n)$, the replacements can only be applied to a depth of $O(1)$. Since we started with $|\cS|=O(n^{1-\delta})$ and each step creates at most $dk$ segments, this terminates after at most $O(n^{1-\delta})$ steps with a state $(B_3',\cS')$ having $O(n^{1-\delta})$ segments, each of cost at most $\alpha n/2$. The total number of holes is $O(n^{1-\delta})=o(n^{1-3\varepsilon})$ and every segment in $\cS'$ lies within distance $\alpha n/2+O(n^{1-\varepsilon})<\alpha n$ of $O$.

Next, we apply \Cref{lem:coverdown-step} once more to every segment of $\cS'$. We arrive at a state $(B_4,\cS_4)$ with the following properties.
\begin{lem} \label{lem:B4}
    $B_4$ is obtained from $B_3$ by poking out a set $A_4$ of $O(n^{1-\delta})$ holes and satisfies the following.
    \begin{enumerate}
        \item $B_4$ is tileable along all but $O(n^{1-\delta})$ lines in $\cL$.
        \item For every line $\ell\in \cL$ along which $B_4$ is not tileable, the following holds.
        \begin{enumerate}
            \item $\ell$ is within distance $\alpha n$ from the center $O$.
            \item $\ell\in \cT(B_3)$.
            \item $\ell\cap B_4$ is a contiguous segment with only one hole, which belongs to some segment of $\cS'$, and is at a distance of at most $\alpha n$ from $O$.
        \end{enumerate}
    \end{enumerate}
\end{lem}
\begin{proof}
    $B_4$ is indeed tileable along all but the lines containing a segment in $\cS_4$, of which there are $O(n^{1-\delta})$ many.

    Each segment $S$ in $\cS_4$ is within distance $\alpha n$ of $O$, so the line $\ell$ containing $S$ is also within distance $\alpha n$ of $O$. From the final applications of \Cref{lem:coverdown-step} to every segment of $\cS'$, we see that $\ell\in \cT(B_3)$ and $\ell\cap B_4$ is a contiguous segment with only one hole, which belongs to some segment of $\cS'$.
\end{proof}

\begin{defn} \label{def+}
    Denote by $\cT_i^+(B_4)$ the union of the set of lines in $\cT_i(B_4)$ (which are pure) and the lines $\ell\in \cL_i$ along which $B_4$ is not tileable.
    Set $\cT^+(B_4)=\bigcup_{i=1}^k \cT_i^+(B_4)$.
\end{defn}

Then we have $\cT(B_4)\subseteq \cT^+(B_4)\subseteq \cT(B_3)$. By \Cref{lem:B4}, $|\cT(B_3)\setminus \cT(B_4)|=O(n^{1-\delta})$ and every line in $\cT^+(B_4)\setminus \cT(B_4)$ is one hole away from being in $\cT(B_4)$, and that hole is within distance $\alpha n$ from $O$.

\section{Tile divisibility} \label{sec:div}

This section implements Step 5 of the proof strategy, which applies our new absorption mechanism.
We start by recalling the definition of tile divisibility, which is the integral relaxation of being tileable.
Henceforth, we will identify a finite subset $B\subset\ZZ^2$ with a Laurent polynomial in $\ZZ[x_1^{\pm},x_2^{\pm}]$ given by $\sum_{b\in B} x^b$.

For a Laurent polynomial $B\in \ZZ[x_1^{\pm},x_2^{\pm}]$ and a line $\ell$, we write $B\cap \ell$ for the polynomial obtained by only keeping the monomials $a_b x^b$ of $B$ for $b \in \ell$.

\begin{defn}
    Given a Laurent polynomial $B\in \ZZ[x_1^{\pm},x_2^{\pm}]$ (or finite subset $B\subset\ZZ^2$), we say that $B$ is \emph{$D_i$-divisible} if $D_i\mid B$. For a line $\ell$ parallel to $v_i$, we say that $B$ is \emph{divisible along $\ell$} if $B\cap \ell$ is $D_i$-divisible.
\end{defn}

For a Laurent polynomial $f=\sum_{u\in\ZZ^2}a_ux^u$, we define its \emph{degree} as $\deg f=\max\setcond{|u_1|+|u_2|}{a_u\neq0}$. Its size is $|f|=\sum_u|a_u|$.

Clearly, if $B$ is tileable by $D_i$, then it is $D_i$-divisible.

At this point, we have $B_4$ which is tileable along all but a small number of lines. In this step, we will poke a small number of holes in $B_4$ to obtain $B_5$, which is divisible along every line (but may not be tileable along a small number of lines). This will be done in $k$ steps. In step $i$, we make $B_4$ divisible along every line parallel to $D_i$, maintaining divisibility for every line parallel to $D_j$ for $j<i$. The main difficulty of this step is ensuring that we only poke a small number of holes, which is resolved via \Cref{thm:poly-decomp}.

\subsection{Two-variable decomposition} \label{subsec:twovar}

We introduce the following useful terminology for decompositions of polynomials into tiny pieces of a specified form.

\begin{defn}
    Let $\cS$ be a class of Laurent polynomials over $\ZZ$. We say that $f$ is \emph{$\cS$-decomposable} if there exist a small number of tiny polynomials $f_1,\ldots,f_m\in\cS$ such that $f=f_1+\cdots+f_m$. We say that the decomposition is \emph{$N$-bounded} if each $f_i$ has degree at most $N$.
\end{defn}

\begin{defn}
    For each $l$, let $\cS_l$ be the set of Laurent polynomials in $\ZZ[x_1^{\pm},x_2^{\pm}]$ which are $D_j$-divisible for $j=1,\ldots,l$.
\end{defn}

We will use the following lemma from \cite{KL26-1}.

\begin{lem}[{\cite[Proposition~5.1]{KL26-1}}] \label{lem:two-var-transference}
    Fix a positive integer $r$ and vectors $u_i=(a_i,b_i)\in \ZZ^2$ with $a_i\neq 0$ for $i=1,\ldots,r$. Fix polynomials $A_i,B_i\in \ZZ[t]$ which are products of (not necessarily distinct) cyclotomic polynomials, with $A_i\mid B_i$ for all $i$, and let $Q\in\ZZ[t]$ be a product of distinct cyclotomic polynomials. In $\ZZ[x^{\pm},y^{\pm}]$, put
    \[
        p=\prod_{i=1}^r A_i(x^{a_i}y^{b_i}),\qquad p^*=\prod_{i=1}^r B_i(x^{a_i}y^{b_i}),\qquad q=Q(y).
    \]
    Then for sufficiently large $n$, if $F\in (p,q)$ is of degree at most $n$, there exists $G\in (p)$ such that $F-G\in (p^*,q)$ and $G$ is a sum of $O(|F|+\log n)$ Laurent polynomials in $(p)$, each of degree at most $n$ and size $O(1)$. All implicit constants may depend on the fixed data, but not on $F$ or $n$.
\end{lem}

\begin{lem} \label{lem:make_divisible}
    Let $1\leq i\leq k$ and let $f\in\cS_{i-1}$. Suppose there exists a small $g\in\ZZ[x_1^{\pm},x_2^{\pm}]$ of degree at most $n$ such that $f-g$ is $D_i$-divisible. Then there exists an $\cS_{i-1}$-decomposable $g^*\in\ZZ[x_1^{\pm},x_2^{\pm}]$ such that $f-g^*\in\cS_i$. Furthermore, the decomposition is $O(n)$-bounded.
\end{lem}

\begin{proof}
    If $i=1$, then we may take $g^*=g$ and decompose it into monomials. Henceforth assume that $i\geq 2$.

    By a change of coordinates, we may assume that $v_i=(0,1)$ and all other $v_j$ have non-zero first coordinates $v_{j1}$. Within this proof we may also reverse $v_j$ so that $v_{j1}>0$. This change of coordinates affects degrees by at most a constant factor, so it is safe to assume without loss of generality that $g$ has degree at most $n$.

    Let
    \[
        p=\prod_{j=1}^{i-1}D_j,
        \qquad q=D_i=1+x_2+\cdots+x_2^{d-1}.
    \]
    The polynomials $D_j$ for distinct directions are pairwise coprime and primitive, so by Gauss's lemma, we have $p\mid f$. Thus $g=f-(f-g)\in(p,q)$.

    Set
    \[
        p^*=p\prod_{j=1}^{i-1}(x^{v_j}-1)
        =\prod_{j=1}^{i-1}(x^{dv_j}-1).
    \]
    Apply \Cref{lem:two-var-transference} to $F=g$, with $r=i-1$, $u_j=v_j$, $x=x_1$, $y=x_2$, and
    \[
        A_j(t)=Q(t)=1+t+\cdots+t^{d-1},\qquad B_j(t)=t^d-1.
    \]
    Then the lemma gives $h\in(p)$ such that $g-h\in(p^*,q)$, together with an $n$-bounded decomposition of $h$ into $O(|g|+\log n)$ tiny polynomials in $(p)$. Since $g$ is small, the number of summands in the decomposition of $h$ is small. Thus $h$ is $\cS_{i-1}$-decomposable and $|h|$ is small.

    Reducing modulo $q$, we can represent every Laurent polynomial $r$ uniquely in the form
    \[
        \tilde r=r_0+x_2r_1+\cdots+x_2^{d-2}r_{d-2},
        \qquad r_j\in\ZZ[x_1^{\pm}].
    \]
    Furthermore, if $r$ is small and has degree $O(n)$, then $\tilde r$ is also small and has degree $O(n)$. Put $r=g-h$. Taking modulo $q$, we have
    \[
        p^*\equiv\overline p
        :=\prod_{j=1}^{i-1}(x_1^{dv_{j1}}-1) \pmod{q},
    \]
    so $(p^*,q)=(\overline p,q)$. Since $r\in (\overline p,q)$, we have $\overline p\mid r_j$ for every $j$.

    By \Cref{thm:poly-decomp} (after shifting $r_j$ by powers of $x_1$ to get an ordinary polynomial), we can write $r_j$ as a small sum of polynomials of the form
    \[
        \pm x_1^a\prod_{l=1}^{i-1}(x_1^{db_lv_{l1}}-1),
        \qquad a\in\ZZ,\quad b_l\geq0,
    \]
    each of degree $O(n)$. Now replace each summand by
    \[
        \pm x_1^a\prod_{l=1}^{i-1}(x^{db_lv_l}-1).
    \]
    This replacement preserves the residue modulo $q$, and the new summand is tiny, belongs to $\cS_{i-1}$, and has degree $O(n)$. We obtain an $O(n)$-bounded $\cS_{i-1}$-decomposable $s_j$ with $q\mid r_j-s_j$.

    Set $s=s_0+x_2s_1+\cdots+x_2^{d-2}s_{d-2}$ and $g^*=h+s$. Then $g^*$ is $\cS_{i-1}$-decomposable with an $O(n)$-bounded decomposition, and
    \[
        f-g^*=f-g+r-s
    \]
    is $D_i$-divisible. Since $f,g^*\in\cS_{i-1}$, we have $f-g^*\in\cS_i$, as required.
\end{proof}

\subsection{Making tile-divisible}

For a polynomial $f\in \ZZ[x_1^{\pm},x_2^{\pm}]$, we denote by $\msupp(f)$ the multiset support of $f$. That is, for every monomial $a_v x^v$ in $f$ with $a_v\neq 0$, we include $v$ in $\msupp(f)$ with multiplicity $|a_v|$.

\begin{defn}[Thin and separated sets]
    Let $A\in\ZZ[x_1^{\pm},x_2^{\pm}]$ be a Laurent polynomial (or a finite (multi)set of points in $\ZZ^2$). Let $v$ be one of the directions $v_1,\ldots,v_k$ and let $\ell$ be a line. Let $C>0$. For polynomials, sizes of intersections count the absolute values of coefficients; unions of multiset supports below add multiplicities.
    \begin{enumerate}
        \item We say that $A$ is \emph{$(C,\ell)$-thin} if $|A\cap \ell|\leq C$. We say that $A$ is \emph{$(C,v)$-thin} if $A$ is $(C,\ell)$-thin for every line $\ell$ parallel to $v$. We say that $A$ is \emph{$C$-thin} if $A$ is $(C,\ell)$-thin for every $\ell\in \cL$.
        \item Define $\supp_v(A)$ to be the set of lines parallel to $v$ that meet (the support of) $A$.
        \item We say that a set of parallel lines is \emph{$C$-separated} if any two distinct lines in the set are at distance at least $C$ from each other. We say that $A$ is \emph{$(C,v)$-projection-separated} if $\supp_v(A)$ is $C$-separated. We say that $A$ is \emph{$C$-projection-separated} if it is $(C,v_j)$-projection-separated for every $j\in[k]$. Note that $A$ is allowed to have many points on the same line, as long as distinct lines are separated.
    \end{enumerate}
\end{defn}

Throughout this subsection, let $\alpha>0$ be a sufficiently small constant relative to $c_{\ref{lem:good-translate}}$ and other necessary constants. To ensure no circular dependencies of the constants, we assert that underlying constants in big-$O$ notation here do not depend on $\alpha$. This $\alpha$ will be used in the cover down step in \Cref{sec:cover-down}.

\begin{lem}[Many good translates] \label{lem:many-translates}
    Let $r\in[k]$ and $C=\Theta(\log^2 n)$. Let $A_1,\ldots,A_m\subset\ZZ^2$, where $m$ is small and each $A_s$ is tiny, lies within distance $O(\alpha n)$ of the origin, and is $(\Omega(\log^2 n),v_j)$-projection-separated for every $j\neq r$. Let $\cM$ be a small set of forbidden lines, each parallel to some $v_j$ with $j\neq r$, and let $\cN\subset\ZZ^2$ be a small set of forbidden points. Then there exist integers $1\leq t_s\leq n/\log n$ such that, writing $A_s'=A_s+t_sdv_r$, the following hold.
    \begin{enumerate}
        \item The sets $A_1',\ldots,A_m'$ are pairwise at distance at least $C$.
        \item $\supp_{v_j}(A_s')\subseteq\cT_j(B_4)$ for every $s$ and every $j\neq r$.
        \item Every point of every $A_s'$ is at distance at least $C$ from every line in $\cM$ and every point in $\cN$.
        \item For distinct $s,s'$ and every $j\neq r$, every line in $\supp_{v_j}(A_s')$ is at distance at least $C$ from every line in $\supp_{v_j}(A_{s'}')$.
    \end{enumerate}
\end{lem}
\begin{proof}
    First consider a single set $A$. For each $j\neq r$ and each line in $\supp_{v_j}(A)$, choose a representative $a_0\in A$ on that line. At each point $a=a_0+qv_j$ on the line, impose the directional requirements $(v_j,-q)$ and $(-v_j,q-1)$, with residues taken modulo $d$. These give a compatible family of requirements, none depending on $v_r$, and ensure that the translated points lie on lines in $\cT_j(B_2)$ for every $j\neq r$.

    Apply \Cref{lem:good-translate} to successive intervals of length $\log^5 n$. This gives $\Omega(n/\log^6 n)$ choices of $1\leq t\leq n/\log n$ for which $A+tdv_r$ satisfies these requirements. Let
    \[
        \cE=\bigcup_{j\neq r}\paren{\cT_j(B_2)\setminus\cT_j(B_4)}.
    \]
    This set is small, since $\cT(B_2)=\cT(B_3)$ and $\cT(B_3)\setminus\cT(B_4)$ is small. Each line in $\cE$ excludes $O(1)$ choices of $t$. Each point in $\cN$ or line in $\cM$ excludes $O(C)$ choices for which some point of $A+tdv_r$ is at distance less than $C$ from it. Thus the total number of exclusions is
    \[
        O\bigl(|\cE|+C(|\cM|+|\cN|)\bigr)=o(n/\log^6 n).
    \]
    Hence $\Omega(n/\log^6 n)$ choices remain.

    In general, we choose $t_1,\ldots,t_m$ successively. When choosing $t_s$, enlarge $\cN$ by $\bigcup_{q<s}A_q'$ and enlarge $\cM$ by
    \[
        \bigcup_{q<s}\ \bigcup_{j\neq r}\supp_{v_j}(A_q').
    \]
    These forbidden sets remain small, so the previous single-set argument applies at every step and gives all the required properties.
\end{proof}

\begin{lem} \label{lem:make-i-good}
    Let $1\leq i\leq k$ and $C=\Theta(\log^2 n)$, and let $g^*$ be $\cS_{i-1}$-decomposable with an $O(\alpha n)$-bounded decomposition. Suppose there exists a small $g\in\ZZ[x_1^{\pm},x_2^{\pm}]$ satisfying the following.
    \begin{enumerate}
        \item $\deg g=O(\alpha n)$.
        \item $g-g^*$ is $D_i$-divisible.
        \item $g$ is $(O(1),v_i)$-thin.
        \item $g$ lies in $\cT_i^+(B_4)$ (recall \Cref{def+}).
    \end{enumerate}
    Then there exists $g'=g_1'+\cdots+g_m'\in\ZZ[x_1^{\pm},x_2^{\pm}]$ satisfying the following.
    \begin{enumerate}
        \item $m$ is small, and each $g_j'\in\cS_{i-1}$ is tiny and has degree $O(\alpha n)$.
        \item $g^*-g'\in\cS_i$. In particular, $g-g'$ is $D_i$-divisible.
        \item Each $g_j'$ is $(C,v_i)$-projection-separated.
        \item $\bigcup_j\msupp(g_j')$ is $(O(1),v_i)$-thin.
        \item Each $g_j'$ lies in $\cT_i^+(B_4)$.
    \end{enumerate}
\end{lem}

\begin{proof}
    Write $g^*=g_1^*+\cdots+g_m^*$, where $m$ is small and each $g_j^*\in\cS_{i-1}$ is tiny and has degree $O(\alpha n)$. Let
    \[
        \cL_0=\supp_{v_i}(g)\cup\bigcup_{j=1}^m\supp_{v_i}(g_j^*).
    \]
    Then $\cL_0$ is small. For each $\ell\in\cL_0$, write
    \[
        g_\ell=g\cap\ell,\qquad g_{j,\ell}^*=g_j^*\cap\ell.
    \]
    Since $g^*-g$ is $D_i$-divisible, so is
    \[
        \sum_{j=1}^m g_{j,\ell}^*-g_\ell.
    \]
    Write $\sum_{j=1}^m g_{j,\ell}^*-g_\ell=\sum_{t\in \ZZ}c_tx^{u+tv_i}$ for some $u\in \ell$. Since this polynomial is $D_i$-divisible, the quantities
    \[\sum_{t\equiv r\pmod d}c_t\]
    are equal for all $r=0,\ldots,d-1$. We can therefore write $g^*-g$ as a small sum of polynomials of the form
    \[
        \pm x^u(1+x^{v_i}+\cdots+x^{(d-1)v_i})
        \quad\text{or}\quad
        \pm x^u(1-x^{dcv_i}),
    \]
    each supported on one line in $\cL_0$, and each of degree $O(\alpha n)$. Let $\cP$ be the collection of these polynomials. For each monomial $bx^u$ in each $g_j^*$, also include $|b|$ copies of each of
    \[
        x^u(1-x^{dv_i})
        \quad\text{and}\quad
        -x^u(1-x^{dv_i})
    \]
    in $\cP$. Thus
    \[
        \sum_{p\in\cP}p=g^*-g,
    \]
    and there are enough monomial occurrences of each sign in $\cP$ to cover all monomial occurrences of the same sign in $g_1^*,\ldots,g_m^*$.

    For each $p\in\cP$, write
    \[
        p=\sum_{u\in P}b_ux^u,
        \qquad b_u\in\{-1,1\},
    \]
    where $P$ lies on a line $\ell_p\in\cL_0$. For each $u\in P$, put
    \[
        q_{p,u}=b_ux^u\prod_{r=1}^{i-1}(1-x^{dc_{p,r}v_r}),
    \]
    where the integers $c_{p,r}=o(n)$ will be chosen below. Let $\cQ$ be the collection of these polynomials, retaining separate copies when necessary.

    We choose the $c_{p,r}$ so that the following hold.
    \begin{enumerate}
        \item Every line in
        \[
            \supp_{v_i}(q_{p,u}-b_ux^u)
        \]
        belongs to $\cT_i^+(B_4)$ and is at distance at least $C$ from every line in $\cL_0$.
        \item For distinct members $p\in\cP$, these sets of lines are disjoint, and their union is $C$-separated. Note that $c_{p,r}$ is independent of the choice of $u\in P$, and that these sets of lines are the same for all $u\in P$.
    \end{enumerate}
    These choices can be made in $i-1$ stages. For each $p$, fix one point $u_p\in P$ and consider the sets
    \[
        S_{p,t}=\setcond{u_p+d\sum_{r\in I}c_{p,r}v_r}{I\subseteq[t]},\quad t=0,1,\ldots,i-1.
    \]
    We maintain that every $S_{p,t}$ is $C$-projection-separated and its points lie on distinct $v_i$-lines, and that the new $v_i$-lines satisfy the two properties above. Initially $S_{p,0}=\{u_p\}$. At stage $r$, let $F_r=\bigcup_{p\in\cP}S_{p,r-1}$ and apply \Cref{lem:many-translates} to the family $(S_{p,r-1})_{p\in\cP}$ in direction $v_r$, with forbidden points $\cN=F_r$ and forbidden lines
    \[
        \cM=\cL_0\cup\bigcup_{l\neq r}\supp_{v_l}(F_r).
    \]
    These sets are small, and no forbidden line is parallel to $v_r$, since $r<i$. This gives all the choices $1\leq c_{p,r}\leq n/\log n$ for this stage. Since
    \[
        S_{p,r}=S_{p,r-1}\cup(S_{p,r-1}+dc_{p,r}v_r),
    \]
    each $S_{p,r}$ remains $C$-projection-separated, and the induction follows.

    Since $c_{p,r}=o(n)$, the $q_{p,u}$ are within $O(\alpha n)+o(n)=O(\alpha n)$ of the origin. Each $q_{p,u}$ is tiny, belongs to $\cS_{i-1}$, and has degree $O(\alpha n)$. Moreover,
    \[
        \sum_{q\in\cQ}q=\sum_{p\in\cP}p\prod_{r=1}^{i-1}(1-x^{dc_{p,r}v_r})\in\cS_i.
    \]
    Each $q_{p,u}$ has exactly one monomial on the lines of $\cL_0$, namely $b_ux^u$.

    We now regroup the polynomials in $\cQ$. For each signed monomial occurrence in each $g_j^*$, select a distinct $q\in\cQ$ whose monomial on $\cL_0$ is that signed monomial. The extra pairs $\pm x^u(1-x^{dv_i})$ added to $\cP$ ensure that this is possible. Let $\cQ_j$ be the collection selected for $g_j^*$, and put
    \[
        h_j=\sum_{q\in\cQ_j}q.
    \]
    Then $|\cQ_j|=O(1)$, and $g_j^*-h_j$ has no support on any line of $\cL_0$.

    Among the unused members of $\cQ$, pair oppositely signed monomials at the same point of $\cL_0$ as far as possible. For every such pair $q,q'$, form the polynomial $-q-q'$. This polynomial also has no support on $\cL_0$.

    The sum of the monomials on $\cL_0$ in the unused members of $\cQ$, before pairing, is $(g^*-g)-g^*=-g$. Consequently, after pairing, if $bx^u$ is a monomial of $g$, exactly $|b|$ unused members remain with monomial $-\operatorname{sgn}(b)x^u$ on $\cL_0$. For each remaining member $q$, consider the polynomial $-q$.

    Relabel all the polynomials considered above
    \[
        g_j^*-h_j,\qquad -q-q',\qquad -q
    \]
    as $g_1',\ldots,g_{m'}'$. Then $m'$ is small. These $g_j'$ satisfy the following, as required.
    \begin{enumerate}
        \item Each $g_j'$ is tiny and belongs to $\cS_{i-1}$.
        \item $g^*-\sum_{j=1}^{m'}g_j'\in\cS_i$, and hence $\sum_{j=1}^{m'}g_j'-g$ is $D_i$-divisible. Indeed, since every member of $\cQ$ is used exactly once,
        \[
            \sum_{j=1}^{m'} g_j'=g^*-\sum_{q\in\cQ}q,
        \]
        and so $\sum_{j=1}^{m'} g_j'-g=(g^*-g)-\sum_{q\in\cQ}q$ is $D_i$-divisible.
        \item Each $g_j'$ is $(C,v_i)$-projection-separated. Indeed, each $g_j'$ has at most one monomial on the lines of $\cL_0$. All its other supporting $v_i$-lines are $C$-separated from each other and from $\cL_0$, by our construction of the $q_{p,u}$.
        \item $\bigcup_j\msupp(g_j')$ is $(O(1),v_i)$-thin. Indeed, for every $\ell\in\cL_0$,
        \[
            \sum_{j=1}^{m'} |g_j'\cap\ell|=|g\cap\ell|=O(1).
        \]
        Every supporting $v_i$-line outside $\cL_0$ comes from some $q_{p,u}$ for only one $p\in\cP$, which contributes only $O(1)$ monomial occurrences.
        \item $\supp(g_j')$ lies in $\cT_i^+(B_4)$.
    \end{enumerate}
\end{proof}

\begin{lem} \label{lem:make-proj-separated}
    Let $1\leq i\leq k$, let $\cM$ be a small set of forbidden lines parallel to $v_i$, and let $N=O(\alpha n)$ and $C=\Theta(\log^2 n)$. Suppose $g\in\cS_{i-1}$ is tiny, has degree at most $N$, lies in $\cT_i^+(B_4)\setminus\cM$, and is $(C,v_i)$-projection-separated. Then there exists $h\in\cS_i$ of size $|h|=O(|g|)$ such that $g-h$ is $C$-projection-separated, lies in $\cT_i^+(B_4)\setminus\cM$, and has degree at most $N+o(n)$.
\end{lem}
\begin{proof}
    The lemma follows from the following inductive step. Let $j\in [k]$. Suppose $g\in\cS_{i-1}$ is tiny, has degree at most $N$, lies in $\cT_i^+(B_4)\setminus\cM$, and is $(C,v_l)$-projection-separated for $l\in[j-1]\cup\{i\}$. Then there exists $h\in\cS_i$ of size $|h|=O(|g|)$ such that $g-h$ is $(C,v_l)$-projection-separated for $l\in[j]\cup\{i\}$, lies in $\cT_i^+(B_4)\setminus\cM$, and has degree at most $N+o(n)$. There is nothing to do when $j=i$. Thus assume $j\neq i$ and that $g$ is already $(C,v_l)$-projection-separated for $l\in[j-1]\cup\{i\}$.

    Let $\ell_1,\ldots,\ell_m$ be the lines of $\supp_{v_j}(g)$, and put $g_t=g\cap\ell_t$. If $j<i$, then $D_j\mid g$ implies $D_j\mid g_t$ for every $t$. Consider
    \[
        f_t=\prod_{\substack{1\leq l\leq i\\l\neq j}}(1-x^{dc_{l,t}v_l}),
        \qquad h=\sum_{t=1}^m g_tf_t,
    \]
    where the integers $c_{l,t}$ will be chosen below. We have $|h|\leq2^i|g|$. If $j<i$, then $D_j\mid g_t$ and $D_l\mid f_t$ for all $l\in [i]\setminus\{j\}$. If $j>i$, then $f_t\in\cS_i$. Hence $h\in\cS_i$ in either case, and
    \[
        g-h=\sum_{t=1}^m g_t(1-f_t).
    \]

    Put $P_t=\supp(g_t)$. Since $P_t$ lies on one $v_j$-line and is $(C,v_i)$-projection-separated, the points of $P_t$ are $C$-separated along its $v_j$-line, and thus $P_t$ is $\Omega(\log^2 n)$-projection-separated in every direction. Start with $U_t=P_t$ and construct the factors of $f_t$ one direction at a time, simultaneously for all $t$.

    When adding the factor in direction $v_r$, where $r\in[i]\setminus\{j\}$, put $U=\bigcup_t U_t$. Apply \Cref{lem:many-translates} to the family $(U_t)_t$ in direction $v_r$, with the forbidden lines
    \[
        \bigcup_{l\neq r}\supp_{v_l}(U),
    \]
    together with $\cM$ if $r\neq i$. This chooses all the $c_{r,t}$ for this stage. Replace each $U_t$ by $U_t\cup(U_t+dc_{r,t}v_r)$.

    Then each $U_t$ remains $\Omega(\log^2 n)$-projection-separated. The same argument shows that $U$ remains $(C,v_l)$-projection-separated for $l\in[j-1]\cup\{i\}$. Moreover, each $U_t$ remains $(C,v_j)$-projection-separated, since this holds initially and $r\neq j$. Thus the newly added $v_j$-lines form a $C$-separated family at distance at least $C$ from $\supp_{v_j}(g)$ and all previously added $v_j$-lines.

    At the end, $g-h$ is supported on $\bigcup_t(U_t\setminus P_t)$, so it has all the required separation and lies in $\cT_i^+(B_4)\setminus\cM$. There are only a bounded number of stages and each $c_{r,t}\leq n/\log n$, so its final degree is at most $N+o(n)$.
\end{proof}

\begin{lem} \label{lem:make-01}
    Let $1\leq i\leq k$, let $\cM$ be a small set of forbidden lines parallel to $v_i$, and let $N=O(\alpha n)$ and $C=\Theta(\log^2 n)$. Suppose $g\in\ZZ[x_1^{\pm},x_2^{\pm}]$ is tiny, has degree at most $N$, is $C$-projection-separated, and lies in $\cT_i^+(B_4)\setminus\cM$. Then there exists a tiny $h\in\cS_k$ such that $g-h$ has only $0$-$1$ coefficients, is $C$-projection-separated, lies in $\cT_i^+(B_4)\setminus\cM$, and has degree at most $N+o(n)$.
\end{lem}
\begin{proof}
    We first describe a polynomial that can be added or subtracted at a prescribed point. Let $A$ be a tiny $C$-projection-separated set lying in $\cT_i^+(B_4)\setminus\cM$ and within distance $O(\alpha n)$ of the origin, and let $u\in A$. We construct
    \[
        F_u=x^u\prod_{p=1}^k
        \left(1+\sum_{t=1}^{d-1}x^{(t+dc_{p,t})v_p}\right),
    \]
    with $1\leq c_{p,t}\leq n/\log n$, so that $F_u$ has only 0-1 coefficients, $\supp(F_u)\cap A=\{u\}$, and $A\cup\supp(F_u)$ remains $C$-projection-separated and lies in $\cT_i^+(B_4)\setminus\cM$.

    Construct the factors for $p=1,\ldots,k$ one at a time, starting with $S=\{u\}$ and maintaining that $A\cup S$ is $C$-projection-separated and lies in $\cT_i^+(B_4)\setminus\cM$. For factor $p$, apply \Cref{lem:many-translates} to the sets $S+tv_p$, $t=1,\ldots,d-1$, in direction $v_p$, with the forbidden lines
    \[
        \bigcup_{l\neq p}\supp_{v_l}(A\cup S),
    \]
    together with $\cM$ if $p\neq i$. This chooses all the $c_{p,t}$ for this factor. The resulting copies $S+(t+dc_{p,t})v_p$ are pairwise disjoint and avoid $A\cup S$. Their union with $A\cup S$ remains $C$-projection-separated. Replace $S$ by the union of $S$ and these copies, then continue with the next factor. In the end, $F_u$ has only 0-1 coefficients, $\supp(F_u)\cap A=\{u\}$, and $A\cup\supp(F_u)$ is $C$-projection-separated and lies in $\cT_i^+(B_4)\setminus\cM$.

    For every $p$, the exponents in its factor run through all residues modulo $d$, so
    \[
        1+\sum_{t=1}^{d-1}x^{(t+dc_{p,t})v_p}
        \equiv 1+x^{v_p}+\cdots+x^{(d-1)v_p}
        \equiv0\pmod{D_p}.
    \]
    Thus $F_u\in\cS_k$, is tiny and is within distance $o(n)$ of $u$.

    Starting with $g$, subtract one such $F_u$ for each positive coefficient occurrence at $u$, using the current support as $A$ at each step. Each subtraction decreases the coefficient at $u$ by one and creates only coefficients $-1$ at new points. After $O(|g|)$ steps, all coefficients are non-positive and the polynomial is still tiny.

    Next, for each negative coefficient occurrence at $u$, add a new $F_u$, again avoiding the current support away from $u$. This increases that negative coefficient by one and creates only coefficients $1$ at new points. It never changes any other existing coefficient. Thus after a bounded number of steps, all coefficients are $0$ or $1$, and the required properties are preserved.
\end{proof}

Combining the above lemmas, we have the following.

\begin{lem} \label{lem:make-div2}
    Let $1\leq i\leq k$ and let $f\in \ZZ[x_1^{\pm},x_2^{\pm}]$ be $D_j$-divisible for $j=1,\ldots,i-1$. Suppose there exists a small $g\in \ZZ[x_1^{\pm},x_2^{\pm}]$ satisfying the following.
    \begin{itemize}
        \item $\deg g=O(\alpha n)$.
        \item $f-g$ is $D_i$-divisible.
        \item $g$ is $(O(1),v_i)$-thin.
        \item $g$ lies in $\cT_i^+(B_4)$.
    \end{itemize}
    Then there exists $g'=g_1'+\cdots+g_m'\in \cS_{i-1}$ satisfying the following.
    \begin{enumerate}
        \item $m$ is small and each $g_j'\in \cS_{i-1}$ is tiny.
        \item $g-g'$ is $D_i$-divisible.
        \item Each $g_j'$ has only 0-1 coefficients and $\deg(g_j')=O(\alpha n)$.
        \item Each $g_j'$ is $\Omega(\log^2 n)$-projection-separated.
        \item $\bigcup_j \msupp(g_j')$ is $(O(1),v_i)$-thin.
        \item Each $g_j'$ lies in $\cT_i^+(B_4)$.
    \end{enumerate}
\end{lem}

\begin{proof}
    Fix $C=\Theta(\log^2 n)$. By \Cref{lem:make_divisible}, there exists an $\cS_{i-1}$-decomposable $g^*$ with an $O(\alpha n)$-bounded decomposition $g^*=g_1^*+\cdots+g_m^*$ such that $g-g^*$ is $D_i$-divisible. Applying \Cref{lem:make-i-good}, we obtain
    \[
        g^{(1)}=g_1^{(1)}+\cdots+g_m^{(1)}
    \]
    with the properties stated there. We reuse the variable $m$ here, which may differ from that in the decomposition of $g^*$, but it is still small.

    Next, apply \Cref{lem:make-proj-separated} successively for $j=1,\ldots,m$, to turn each $g_j^{(1)}$ into $g_j^{(2)}$ which is additionally $C$-projection-separated. More precisely, at step $j$, use the small forbidden set
    \[
        \cM_j^{(1)}=\paren{\bigcup_{r=1}^m\supp_{v_i}(g_r^{(1)})\cup\bigcup_{r<j}\supp_{v_i}(g_r^{(2)})}\setminus\supp_{v_i}(g_j^{(1)}).
    \]
    \Cref{lem:make-proj-separated} gives a tiny $h_j^{(1)}\in\cS_i$ such that
    \[
        g_j^{(2)}=g_j^{(1)}-h_j^{(1)}
    \]
    is $C$-projection-separated and lies in $\cT_i^+(B_4)\setminus\cM_j^{(1)}$. These $g_j^{(2)}$ satisfy the following.
    \begin{enumerate}
        \item Each $g_j^{(2)}$ is tiny and belongs to $\cS_{i-1}$.
        \item $\sum_{j=1}^m g_j^{(2)}-g$ is $D_i$-divisible. Indeed, both $\sum_{j=1}^m g_j^{(1)}-g$ and each $g_j^{(2)}-g_j^{(1)}=-h_j^{(1)}$ are $D_i$-divisible.
        \item Each $g_j^{(2)}$ is $C$-projection-separated.
        \item $\bigcup_j\msupp(g_j^{(2)})$ is $(O(1),v_i)$-thin. Indeed, $\bigcup_j\msupp(g_j^{(1)})$ is $(O(1),v_i)$-thin. Any line parallel to $v_i$ meeting some $g_r^{(1)}$ can meet $g_j^{(2)}$ only if it already meets $g_j^{(1)}$, so there are only $O(1)$ such indices $j$. Any line meeting none of the $g_r^{(1)}$ can meet at most one $g_j^{(2)}$.
        \item $\supp(g_j^{(2)})$ lies in $\cT_i^+(B_4)$.
    \end{enumerate}

    Finally, apply \Cref{lem:make-01} successively, using the small forbidden set
    \[
        \cM_j^{(2)}=\paren{\bigcup_{r=1}^m\supp_{v_i}(g_r^{(2)})\cup\bigcup_{r<j}\supp_{v_i}(g_r')}\setminus\supp_{v_i}(g_j^{(2)}).
    \]
    We obtain tiny $h_j^{(2)}\in\cS_k$ such that
    \[
        g_j'=g_j^{(2)}-h_j^{(2)}
    \]
    has the same properties as $g_j^{(2)}$ listed above, in addition to having only 0-1 coefficients. Setting $g'=\sum_jg_j'$ gives all the required properties.
\end{proof}

\begin{lem} \label{lem:make-i-div}
    Let $1\leq i\leq k$ and let $f\in \cS_{i-1}$. Suppose there exists a small $g\in \ZZ[x_1^{\pm},x_2^{\pm}]$ such that $f-g$ is $D_i$-divisible and $g$ is $(O(1),v_i)$-thin, lies in $\cT_i^+(B_4)$, and has degree $O(\alpha n)$. Let $A\subset \ZZ^2$ be a small set. Then there exists a small $g' \in \cS_{i-1}$ satisfying the following.
    \begin{enumerate}
        \item $g-g'$ is $D_i$-divisible.
        \item $g'$ has only 0-1 coefficients.
        \item $\deg(g')=O(\alpha n)$.
        \item $\supp(g')$ is disjoint from $A$.
        \item $\supp(g')$ is $O(1)$-thin.
        \item $\supp(g')$ lies in $\cT_j^+(B_4)$ for all $j$.
    \end{enumerate}
\end{lem}
\begin{proof}
    By \Cref{lem:make-div2}, we obtain $g''=g_1''+\cdots+g_m''\in \cS_{i-1}$. We will set $g_j'$ to be a translate $g_j'=x^{c_j dv_i} g_j''$ for some positive integer $c_j$, and take $g'=\sum_j g_j'$. For any choice of $c_j\leq n/\log n$, the following hold.
    \begin{enumerate}
        \item $g_j'$ is tiny and in $\cS_{i-1}$.
        \item $g_j'-g_j''$ is $D_i$-divisible, so $g-g'$ is $D_i$-divisible.
        \item $g_j'$ has only 0-1 coefficients.
        \item $\deg(g_j')=O(\alpha n)$.
        \item $\supp(g')$ is $(O(1), v_i)$-thin.
        \item $\supp(g')$ lies in $\cT_i^+(B_4)$.
    \end{enumerate}

    Apply \Cref{lem:many-translates} to the family $(\supp(g_j''))_{j=1}^m$ in direction $v_i$, with $C=\Theta(\log^2 n)$, $\cM=\emptyset$, and $\cN=A$. This chooses all the $c_j$ simultaneously. The supports of the $g_j'$ are pairwise disjoint and avoid $A$, so $g'$ has only $0$-$1$ coefficients and avoids $A$. For every $l\neq i$, these supports lie on lines in $\cT_l(B_4)$, and their sets of supporting $v_l$-lines are pairwise disjoint. Since each $g_j'$ is tiny, their union is $(O(1),v_l)$-thin. Together with the properties above in direction $v_i$, this proves the lemma.
\end{proof}

\begin{lem}[Construction of $B_5$] \label{lem:B5}
    There exists a $B_5$ obtained by poking a small set $A_5$ of holes in $B_4$ such that the following holds.
    \begin{enumerate}
        \item $A_5$ lies in $\cT_i^+(B_4)$ for all $i$.
        \item $A_5$ is $O(1)$-thin.
        \item $A_5$ is within distance $O(\alpha n)$ from the origin.
        \item $B_5$ is $D_i$-divisible for all $i=1,\ldots,k$.
    \end{enumerate}
\end{lem}

\begin{proof}
    Let $H_i$ be the set of holes on the lines in $\cT_i^+(B_4)\setminus\cT_i(B_4)$. By \Cref{lem:B4}, each such line has exactly one hole. This hole is within distance $\alpha n$ from the origin, and $|H_i|=O(n^{1-\delta})$. All other lines parallel to $v_i$ are already tileable, so $B_4+H_i$ is $D_i$-divisible.

    We will apply \Cref{lem:make-i-div} $k$ times. Set $B_4^0=B_4$. Suppose after $i-1$ steps we have $B_4^{i-1}=B_4\setminus(F_1\cup\cdots\cup F_{i-1})$, where the sets $F_r$ are small and satisfy the following.
    \begin{enumerate}
        \item $B_4^{i-1}$ is $D_j$-divisible for $j=1,\ldots,i-1$.
        \item $F_1,\ldots,F_{i-1}$ are disjoint and lie within distance $O(\alpha n)$ from the origin.
        \item $F_1\cup\cdots\cup F_{i-1}$ is $O(1)$-thin and lies in $\cT_j^+(B_4)$ for all $j$.
    \end{enumerate}
    Let $g=-H_i-F_1-\cdots-F_{i-1}$. Then $B_4^{i-1}-g=B_4+H_i$ is $D_i$-divisible. Since $H_i$ is $(1,v_i)$-thin, $g$ is $(O(1),v_i)$-thin.

    Apply \Cref{lem:make-i-div} with $f=B_4^{i-1}$ and $A=F_1\cup\cdots\cup F_{i-1}\cup\bigcup_{j=1}^kH_j$. Its hypotheses hold because $g$ is small, has degree $O(\alpha n)$, and lies in $\cT_i^+(B_4)$. We obtain a small $g'\in\ZZ[x_1^{\pm},x_2^{\pm}]$ such that the following hold.
    \begin{enumerate}
        \item $g'$ is $D_j$-divisible for $j=1,\ldots,i-1$.
        \item $g-g'$ is $D_i$-divisible.
        \item $g'$ has only $0$-$1$ coefficients and degree $O(\alpha n)$.
        \item Let $F_i=\supp(g')$. Then $F_i$ is disjoint from $F_1\cup\cdots\cup F_{i-1}$ and $\bigcup_j H_j$.
        \item $F_i$ is $O(1)$-thin.
        \item $F_i$ lies in $\cT_j^+(B_4)$ for all $j$.
    \end{enumerate}
    The points of $F_i$ lie in the central region of $B_1$. There, a point on a line in $\cT_j^+(B_4)$ belongs to $B_4$ unless it is in $H_j$. Thus avoidance of $\bigcup_jH_j$ and of the previous sets $F_r$ gives $F_i\subset B_4^{i-1}$. Set $B_4^i=B_4^{i-1}\setminus F_i$. By the properties above, we have the following.
    \begin{enumerate}
        \item $B_4^i$ is $D_j$-divisible for $j=1,\ldots,i-1$ since both $B_4^{i-1}$ and $g'$ are.
        \item $B_4^i=(B_4^{i-1}-g)+(g-g')$ is $D_i$-divisible. Hence $B_4^i$ is $D_j$-divisible for $j=1,\ldots,i$.
        \item $F_i$ is disjoint from $F_1\cup\cdots\cup F_{i-1}$ and is within distance $O(\alpha n)$ from the origin.
        \item $F_1\cup\cdots \cup F_i$ is $O(1)$-thin.
        \item $F_1\cup\cdots\cup F_i$ lies in $\cT_j^+(B_4)$ for all $j$.
    \end{enumerate}
    This gives the induction step. Setting $B_5=B_4^k$ and $A_5=F_1\cup\cdots\cup F_k$ gives the desired result.
\end{proof}

\section{From tile divisibility to tileability} \label{sec:tile}

It remains to implement Step 6 of the proof strategy.
Recall that from \Cref{lem:B5}, we have a set $B_5=B_4\setminus A_5$ satisfying the following.
\begin{enumerate}
    \item $A_5$ lies in $\cT_i^+(B_4)$ for all $i$.
    \item $A_5$ is small and $O(1)$-thin.
    \item $A_5$ is within distance $O(\alpha n)$ from the origin.
    \item $B_5$ is $D_i$-divisible for all $i=1,\ldots,k$.
\end{enumerate}
Let $\cN_i$ be the set of lines parallel to $v_i$ intersecting $A_5$. Then $B_5$ is tileable along every line in $\cL_i\setminus \cN_i$. Every line $\ell\in \cN_i$ contains only $O(1)$ many points of $A_5$, since $A_5$ is $O(1)$-thin. Furthermore, $\ell$ contains at most one other hole in $B_5$ which does not belong to $A_5$, and is within distance $O(\alpha n)$ from the origin.

To make every line tileable, the idea is to design a gadget that allows us to move the holes in $A_5$ along the lines in steps of multiples of $d$. Since every line is tile-divisible, there exists a permutation of the ordering of the holes that makes the line tileable. Indeed, if $B_5$ is tile-divisible along some $\ell\in \cT_i(B_3)$ where $B_5\cap \ell$ has holes $a_1,\ldots,a_m$, then the set of holes $\{a_1,\ldots,a_m\}$ is itself tile-divisible. Thus, the collection of $\tau$-distances $\tau_{v_i}(a_j)$ consists of exactly $m/d$ copies of each residue modulo $d$. By moving the holes in steps of multiples of $d$, we can achieve any ordering of the holes along $\ell$, so we can find an ordering that makes $B_5$ tileable along $\ell$.

We now prescribe a movement plan for each hole in $A_5$ and each direction. Let $1\leq i\leq k$ and $\ell\in \cN_i$ with holes $a_1,\ldots,a_m$, all of which belong to $A_5$ except possibly one of them. Since $A_5$ is $O(1)$-thin, we have $m=O(1)$.

First suppose $a_j\in A_5$ for all $j$. From the above discussion, we can move the holes $a_1,\ldots,a_m$ in multiples of $d$ so that the resulting ordering of the holes along $\ell$ makes $\ell$ tileable. Hence, we can define disjoint segments $S_{a_1,i},\ldots,S_{a_m,i}$ on $\ell$, each with its own set of permissible cells, satisfying the following properties.
\begin{itemize}
    \item A cell $p\in S_{a_j,i}$ is permissible if and only if $p-a_j\in d\ZZ v_i$.
    \item Each $S_{a_j,i}$ has length $\Omega(n/\log n)$ and does not contain $a_j$.
    \item For any choice of permissible cells $a_j'\in S_{a_j,i}$, $B_5\cup \set{a_1,\ldots,a_m}\setminus \set{a_1',\ldots,a_m'}$ is tileable along $\ell$.
\end{itemize}

If there is one $a_s$ that does not belong to $A_5$, leave it fixed and prescribe segments $S_{a_j,i}$ for all other $a_j$ as above. In both cases, since there are only $O(1)$ original holes, we can choose pairwise disjoint segments of length $\Theta(n/\log n)$ avoiding all of them. All these segments lie within distance $O(\alpha n)$ of the origin. This gives the following lemma.

\begin{lem} \label{lem:move-plan}
    For each $a\in A_5$ and $1\leq i\leq k$, there exists a segment $S_{a,i}$ of length $\Omega(n/\log n)$, lying within distance $O(\alpha n)$ of the origin, on the line through $a$ parallel to $v_i$. Declare a cell $p\in S_{a,i}$ permissible if and only if $p-a\in d\ZZ v_i$. These segments have the following property. For any line $\ell$ parallel to $v_i$, suppose $\ell\cap A_5=\{a_1,\ldots,a_m\}$. Then the segments $S_{a_j,i}$ are disjoint and avoid all the holes on $\ell$, and for any choice of permissible cells $a_j'\in S_{a_j,i}$,
    \[
        B_5\cup \set{a_1,\ldots,a_m}\setminus \set{a_1',\ldots,a_m'}
    \]
    is tileable along $\ell$.
\end{lem}

Since $S_{a,i}$ does not contain $a$, the segment lies in the direction $\epsilon_{i} v_i$ from $a$ for some $\epsilon_{i}\in \{\pm 1\}$. The sign $\epsilon_i$ depends on $a$, but we suppress this dependence for ease of notation.

For each $a\in A_5$ and $1\leq i\leq k$, let $u_{a,i}$ be a permissible cell in $S_{a,i}$ and define the polynomial
\[
    g_a=x^a\prod_{i=1}^k (1-x^{u_{a,i}-a}).
\]
The monomials of $g_a$ are of the form $(-1)^{|T|}x^{a(T)}$ for some subset $T\subseteq [k]$, where $a(T)=a+\sum_{i\in T} (u_{a,i}-a)$. Let $\triangle$ denote the symmetric difference. All these points lie within distance $O(\alpha n)$ of the origin, with underlying constants independent of $\alpha$. Take $\alpha$ sufficiently small that this region and its $O(n/\log n)$-neighborhood lie in $B_1$ and in the ball of radius $c_{\ref{lem:good-translate}}n$, so that \Cref{lem:good-translate} is applicable.

Some of the monomials of $g_a$ across all $a$ trivially lie on the same line. For instance, for a line $\ell\in \cN_i$ and $a,a'\in \ell\cap A_5$, the pairs of monomials $x^a-x^{u_{a,i}}$ of $g_a$ and $x^{a'}-x^{u_{a',i}}$ of $g_{a'}$ lie on the same line $\ell$. For $a\in A_5$ and $T\subseteq [k]$, another such collision occurs because the monomials $(-1)^{|T|}x^{a(T)}$ and $(-1)^{|T|+1}x^{a(T\triangle \{i\})}$ of $g_a$ lie on the same line parallel to $v_i$. We say that two distinct terms, labelled by $(a,T)$, have a \emph{collision in direction $v_i$} if their points $a(T)$ lie on the same line parallel to $v_i$. The two types of collisions described above are called \emph{trivial collisions} in that direction. We can choose the $u_{a,i}$ sequentially to avoid all non-trivial collisions. The following lemma says that the $u_{a,i}$ can be chosen to satisfy additional properties.
\begin{lem} \label{lem:good-hypercube}
    There exists a choice of $u_{a,i}$ such that $\supp_{v_j}(g_a)\subseteq\cT_j(B_5)\cup\{a+\RR v_j\}$ for every $a$ and $j$, and the family $(g_a)_{a\in A_5}$ has no non-trivial collision. Furthermore, for each monomial of the form $x^{a(T)}$ with $T\neq \emptyset$ and $|T|$ even, and every $i$,
    \begin{itemize}
        \item if $i\in T$, then $\tau_{-\epsilon_{i} v_i}(a(T))=0$;
        \item if $i\notin T$, then $\tau_{\epsilon_{i} v_i}(a(T))=0$.
    \end{itemize}
    In other words, the $\tau$ distance from $a(T)$ towards $a(T\triangle \{i\})$ is 0 mod $d$.
\end{lem}
\begin{proof}
    Fix $a\in A_5$. For $|T|$ even, let $P(T,i)$ be the assertion that the $\tau$-distance from $a(T)$ towards $a(T\triangle \{i\})$ is 0 mod $d$, and that $a(T)$ lies on $\cT_i(B_3)$. Then $P(T,i)$ is equivalent to a requirement on $a(T)$ (see \Cref{lem:good-translate}). Since $a(T\triangle \{i\})-a(T)\in d\ZZ v_i$, $P(T,i)$ is also equivalent to a similar requirement on $a(T\triangle \{i\})$, which we call $P(T\triangle \{i\},i)$.

    We will find $u_{a,i}$ for $i=1,\ldots,k$ sequentially. Let $g_i$ be the polynomial $x^a\prod_{j=1}^i (1-x^{u_{a,j}-a})$, with $g_0=x^a$. We claim that we can find $u_{a,i}$ such that the following holds for all $i$.
    \begin{enumerate}
        \item $g_i$ is $\Omega(\log^2 n)$-projection-separated.
        \item There is no non-trivial collision within $g_i$ or between $g_i$ and any previously constructed $g_{a'}$.
        \item $\supp_{v_j}(g_i)\subseteq\cT_j(B_5)\cup\{a+\RR v_j\}$ for all $j$.
        \item $P(T,j)$ holds for all $\emptyset\neq T\subseteq [i]$ and all $j$, except for $T=\{j\}$.
    \end{enumerate}

    Suppose we have found $u_{a,1},\ldots,u_{a,i-1}$, and we are finding $u_{a,i}$. For $T\subseteq [i]$ containing $i$ and $j\neq i$, the condition $P(T,j)$ is equivalent to a requirement in the directions $\pm v_j$ on the translate $x^{u_{a,i}-a}g_{i-1}$. It can be checked that these requirements are compatible. Applying \Cref{lem:good-translate} to successive intervals of length $O(\log^5 n)$ in the permissible cells of $S_{a,i}$ gives $\Omega(n/\log^6 n)$ choices of $u_{a,i}$ satisfying all these requirements. Since $P(T,i)\iff P(T\triangle \{i\},i)$, the requirements in direction $v_i$ are inherited, except for $T=\emptyset,\{i\}$. Thus $P(T,j)$ holds for all $\emptyset\neq T\subseteq [i]$ and all $j$, except for $T=\{j\}$.

    Of these $\Omega(n/\log^6 n)$ choices of $u_{a,i}$, only a small number of them are bad. Indeed, the bad choices are those satisfying one of the following.
    \begin{enumerate}
        \item $g_i$ is not $\Omega(\log^2 n)$-projection-separated or has a non-trivial collision with itself. There are only $O(|g_{i-1}|^2\log^2 n)$ many such bad choices.
        \item There is a non-trivial collision between $g_i$ and some previously constructed $g_{a'}$. In each direction $v_j$ with $j\neq i$, every pair of terms excludes at most one choice, and there are only a small number of pairs.
        \item $\supp_{v_j}(g_i)\not\subseteq\cT_j(B_5)\cup\{a+\RR v_j\}$ for some $j$. This cannot happen for $j=i$. For $j\neq i$, the translated terms already lie in $\cT_j(B_3)$, so failure requires meeting one of the small number of lines in $\cT_j(B_3)\setminus\cT_j(B_5)$. Each such line excludes only $O(1)$ choices.
    \end{enumerate}

    The total number of exclusions is small, hence $o(n/\log^6 n)$. Thus, we can find $u_{a,i}$ satisfying all the required properties.

\end{proof}

Consider a monomial of $g_a$ with coefficient $+1$, other than $x^a$. Such a monomial is of the form $x^{a(T)}$ for some non-empty $T\subseteq [k]$ of even size. For each $1\leq i\leq k$, let $\ell_i$ be the line through $a(T)$ parallel to $v_i$. Then the only other monomial of $g_a$ that is on $\ell_i$ is $x^{a(T\triangle \{i\})}$, which has coefficient $-1$. The point $a(T\triangle \{i\})$ is in the direction of $\epsilon_i' v_i$ from $a(T)$, where
\[\epsilon_i'=\begin{cases}
        -\epsilon_{i} & \text{if } i\in T     \\
        \epsilon_{i}  & \text{if } i\notin T.
    \end{cases}\]

For each $a\in A_5$ and non-empty even-sized $T\subseteq[k]$, define the \emph{corner gadget} $g_{a,T}$ to be the polynomial
\[
    g_{a,T}=x^{a(T)}\prod_{i=1}^k (1+x^{-\epsilon_i' c_{i,1}(T) v_i}+\cdots+x^{-\epsilon_i' c_{i,d-1}(T) v_i}),
\]
for some increasing positive integers $c_{i,j}(T)\equiv j\pmod{d}$ for $j=1,\ldots,d-1$ to be chosen later, with $c_{i,j}(T)=O(n/\log n)$. The parameters $c_{i,j}(T)$ also depend on $a$, and $\epsilon_i'$ depends on both $a$ and $T$, but we suppress these dependencies for ease of notation. Put $c_{i,0}(T)=0$, and label the monomials of $g_{a,T}$ by tuples $(j_1,\ldots,j_k)\in\{0,1,\ldots,d-1\}^k$, with corresponding exponent
\[
    a(T)-\sum_{i=1}^k\epsilon_i' c_{i,j_i}(T)v_i.
\]
Just as for $g_a$, the terms of $g_{a,T}$ can have unavoidable collisions with each other and with $g_a$. More precisely, we have the following additional \emph{trivial collisions} in direction $v_i$ (other than the ones already defined for $g_a$).
\begin{enumerate}
    \item The $d$ monomials of $g_{a,T}$ with labels $(j_1,\ldots,j_k)$, for each fixed choice of $j_r$ for $r\neq i$, and $j_i=0,1,\ldots,d-1$, all lie on the same line parallel to $v_i$.
    \item The $d$ monomials of $g_{a,T}$ with $j_r=0$ for every $r\neq i$, and the two monomials of $g_a$ at $a(T)$ and $a(T\triangle\{i\})$, all lie on the same line through $a(T)$ parallel to $v_i$.
    \item The $x^{a(T)}$ monomials of $g_{a,T}$ and $g_a$ coincide.
\end{enumerate}
All other potential collisions are non-trivial. It is not hard to see that there are choices of the $c_{i,j}(T)$ that avoid all non-trivial collisions (say by choosing randomly).

\begin{figure}
    \centering
    \includegraphics[scale=1.0]{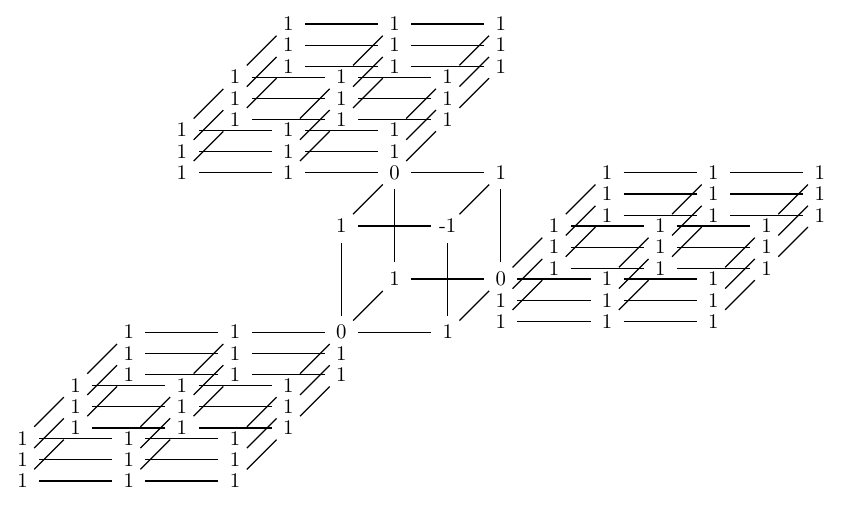}
    \caption{Coefficients of $h_a$, for $k=d=3$. The center cube is $-g_a$, with coefficient $-1$ on the monomial $x^a$. The outer cubes are the $g_{a,T}$, and always point away from the center cube.}
    \label{fig:ha}
\end{figure}

Let $h_a=-g_a+\sum_T g_{a,T}$, where the sum is over all non-empty even-sized subsets of $[k]$. If there are no non-trivial collisions, then all nonzero coefficients of $h_a$ are $+1$, except for $x^a$ which has coefficient $-1$. See \Cref{fig:ha} for an example illustration of $h_a$.

\begin{lem} \label{lem:corner-gadget}
    There exists a choice of $c_{i,j}(T)$ such that the following hold.
    \begin{enumerate}
        \item $g_{a,T}$ has support in $\cT_j(B_5)$ for all $j$.
        \item For any line $\ell\in \cL_j$ intersecting $g_{a,T}$, $B_2-g_{a,T}$ is tileable along $\ell$.
        \item There are no non-trivial collisions. In particular, $h_a$ has no non-trivial collision with itself or any other $h_{a'}$.
    \end{enumerate}
\end{lem}
\begin{proof}
    The monomials of $g_{a,T}$ are on the points of the form
    \[
        a(T)-\sum_{i=1}^k \epsilon_i' c_{i,j_i}(T) v_i
    \]
    for some $j_i\in \{0,1,\ldots,d-1\}$, where $c_{i,0}(T)=0$. To ensure that $B_2-g_{a,T}$ is tileable along each $\ell$ intersecting $g_{a,T}$, we require that
    \begin{equation} \label{eq:corner-req}
        \tau_{\epsilon_l' v_l}\paren{a(T)-\sum_{i=1}^k \epsilon_i' c_{i,j_i}(T) v_i}\equiv 0\pmod{d}
    \end{equation}
    for all $l$ such that $j_l=0$.

    Define the partially constructed polynomials
    \[
        g_{a,T,i}=x^{a(T)}\prod_{j=1}^i (1+x^{-\epsilon_j' c_{j,1} v_j}+\cdots+x^{-\epsilon_j' c_{j,d-1} v_j})
    \]
    for $i=0,\ldots,k$, where we suppress the dependence on $T$ by writing $c_{i,j}=c_{i,j}(T)$ for the remainder of this proof. Then $g_{a,T,k}=g_{a,T}$. Observe that $g_{a,T,i+1}$ consists of $g_{a,T,i}$ and $d-1$ other translates of $g_{a,T,i}$, by the vectors $-\epsilon_{i+1}' c_{i+1,j}v_{i+1}$ for $j=1,\ldots,d-1$. We will find $c_{i,j}$ such that the partially constructed polynomials $g_{a,T,i}$ satisfy the following properties for all $i$.
    \begin{enumerate}
        \item $g_{a,T,i}$ has support in $\cT_j(B_5)$ for all $j$.
        \item \eqref{eq:corner-req} holds at every point for each $l$ with $j_l=0$, taking $j_l=0$ for $l>i$.
        \item There is no non-trivial collision.
        \item $g_{a,T,i}$ is $\Omega(\log^2 n)$-projection-separated.
    \end{enumerate}
    For the base case, $g_{a,T,0}=x^{a(T)}$ satisfies all the properties. Now suppose we have $g_{a,T,i}$. We will choose $c_{i+1,1}<c_{i+1,2}<\cdots<c_{i+1,d-1}$ sequentially, with $c_{i+1,j}\equiv j\pmod{d}$ for $1\leq j\leq d-1$, $c_{i+1,1}=O(n/\log n)$, and $c_{i+1,j+1}-c_{i+1,j}=O(n/\log n)$ for $1\leq j<d-1$.

    This is another application of \Cref{lem:good-translate}. Finding a suitable $c_{i+1,j}$ is equivalent to finding a suitable good translate of $g_{a,T,i}$ along the direction $-\epsilon_{i+1}'v_{i+1}$. Having seen many similar applications, we omit the details here.
\end{proof}

Fix a choice of $c_{i,j}(T)$ as above. This gives us $h_a$ for each $a\in A_5$. The key property of $h_a$ is the following.

\begin{lem} \label{lem:ha-prop}
    Let $\ell$ be a line parallel to $v_i$ which intersects $h_a$ for some $a\in A_5$. If $\ell$ contains $a$, then $h_a\cap \ell$ is just $-x^a+x^{u_{a,i}}$. Otherwise, $\ell$ does not intersect any other $h_{a'}$ for $a'\neq a$, and $h_a\cap \ell$ is such that $B_5-h_a$ is tileable along $\ell$.
\end{lem}
\begin{proof}
    If $\ell$ contains $a$, then from the lack of non-trivial collisions, the only other monomial of $h_a$ on $\ell$ is $x^{u_{a,i}}$, coming from $-g_a$. If $\ell$ does not contain $a$, we consider separately the cases where $\ell$ intersects $g_a$ or not.

    If $\ell$ intersects $g_a$, then $\ell\cap g_a$ is of the form $x^{a(T)}-x^{a(T\triangle \{i\})}$ for some non-empty even-sized $T\subseteq [k]$. Then $\ell\cap g_{a,T}$ is of the form
    \[x^{a(T)}+x^{a(T)-\epsilon_i' c_{i,1}(T) v_i}+\cdots+x^{a(T)-\epsilon_i' c_{i,d-1}(T) v_i}.\]
    Thus, from the lack of non-trivial collisions, $\ell\cap h_a=\ell\cap (-g_a)+\ell\cap g_{a,T}$ is of the form
    \[x^{a(T\triangle \{i\})}+x^{a(T)-\epsilon_i' c_{i,1}(T) v_i}+\cdots+x^{a(T)-\epsilon_i' c_{i,d-1}(T) v_i}.\]
    Since $a(T\triangle \{i\})$ is in the direction of $\epsilon_i' v_i$ from $a(T)$, the monomials are listed in order of the direction $-\epsilon_i' v_i$ along $\ell$. By \Cref{lem:good-hypercube}, we have $\tau_{\epsilon_i' v_i}(a(T))=0$, so $\tau_{\epsilon_i' v_i}(a(T\triangle \{i\}))\equiv 0\pmod{d}$. Since $c_{i,j}(T)\equiv j\pmod{d}$, we have $\tau_{\epsilon_i' v_i}(a(T)-\epsilon_i' c_{i,j}(T) v_i)\equiv j\pmod{d}$ for all $j$. Thus, $B_5-h_a$ is tileable along $\ell$.

    If $\ell$ does not intersect $g_a$, then it intersects a unique $g_{a,T}$ for some $T$. By \Cref{lem:corner-gadget} and the lack of non-trivial collisions, $(B_5-g_{a,T})\cap\ell=(B_2-g_{a,T})\cap\ell$ is tileable along $\ell$.
\end{proof}

In essence, the above lemma says that by replacing $B_5$ with $B_5-h_a$, we have in effect moved the hole $a$ to $u_{a,i}$ along the line parallel to $v_i$, while preserving tileability along all other lines. Finally, set $B_6=B_5-\sum_{a\in A_5} h_a$. Because there are no non-trivial collisions, $B_6$ has only 0-1 coefficients, so is an actual subset of $B_2$. This is our construction for \Cref{thm:main-ap}.

\begin{proof}[Proof of \Cref{thm:main-ap}]
    We claim that $B_6$ works. Recall that $B_5$ is tileable along every line not intersecting $A_5$. Thus the only possible line for which $B_6$ is not tileable is a line intersecting $A_5$ or some $h_a$. Consider such a line $\ell$, which is, say, parallel to $v_1$.

    If $\ell$ does not intersect $A_5$, then it intersects exactly one $h_a$. By \Cref{lem:ha-prop}, $B_6\cap \ell=(B_5-h_a)\cap \ell$ is tileable.

    Otherwise $\ell$ intersects $A_5$. Write $A_5\cap \ell=\{a_1,\ldots,a_m\}$. Recall that $B_5\cap \ell$ either has holes $A_5\cap\ell$, or has one other hole, denoted by $a_0$, within distance $O(\alpha n)$ of the origin. If $a_0$ does not exist, then by \Cref{lem:ha-prop}, the holes of $B_6\cap \ell$ are $u_{a_1,1},\ldots,u_{a_m,1}$. These holes lie in the permissible cells of the segments $S_{a_1,1},\ldots,S_{a_m,1}$, so by construction, $B_6$ is tileable along $\ell$. If $a_0$ exists, then the holes of $B_6\cap \ell$ are $a_0,u_{a_1,1},\ldots,u_{a_m,1}$. Again, by the construction of the segments $S_{a_j,1}$, $B_6$ is tileable along $\ell$.

\end{proof}

\section{Concluding remarks}  \label{sec:conclude}

We conclude by briefly discussing some natural variants of \Cref{prob:dominoes}.
\begin{itemize}
    \item If $B$ is allowed to be infinite, one can simply take $B=\ZZ^2$.
    \item For tiles other than dominoes and APs, the answer is generally negative. For instance, consider the following two tiles:
    \begin{center}
        \begin{tikzpicture}[x=0.45cm, y=0.45cm, line width=0.6pt]
            \draw (0,0) rectangle (2,1);
            \draw (1,0) -- (1,1);
            \node at (1,-0.5) {$D_1$};
            \begin{scope}[shift={(6,0)}]
                \draw (0,0) -- (2,0) -- (2,1) -- (1,1)
                -- (1,2) -- (0,2) -- cycle;
                \draw (1,0) -- (1,1) -- (0,1);
                \node at (1,-0.5) {$D_2$};
            \end{scope}
        \end{tikzpicture}
    \end{center}
    To see that no finite $B$ is simultaneously tileable by both, let $X$ be the rightmost of the highest cells of $B$.
    To be tileable by $D_1$, $B$ must contain the cell $Y$ immediately to the left of $X$.
    For $B$ to be tileable by $D_2$, $X$ must be covered by the top cell of a copy of $D_2$. But now, no copy of $D_2$ can cover $Y$ without overlapping.
    \item In one dimension (dominoes belonging to $\ZZ$), the answer is easily positive by considering a long segment. The problem is significantly more complex in higher dimensions, but we conjecture that the answer is still positive.

\end{itemize}
\begin{conj} \label{conj:higher-dim}
    Let $d\geq 3$ and $D_1,\ldots,D_k\subset\ZZ^d$ be APs. Then there exists a
    finite non-empty $B\subset\ZZ^d$ that is tileable by every $D_i$.
\end{conj}

\bibliographystyle{plain}
\bibliography{main}

@inproceedings{keevash2019hypergraph,
  title        = {Hypergraph matchings and designs},
  author       = {Keevash, Peter},
  booktitle    = {Proceedings of the International Congress of Mathematicians (ICM 2018)},
  volume       = {3},
  year         = {2019},
  organization = {World Scientific Publishing}
}

@unpublished{KL26-1,
  author = {Keevash, Peter and Lim, Jeck},
  title  = {Efficient reconfiguration, polynomial decompositions and approximate counting},
  year   = {2026},
  note   = {preprint}
}

@article{EGJ20,
  title   = {Pseudorandom hypergraph matchings},
  volume  = {29},
  doi     = {10.1017/S0963548320000280},
  number  = {6},
  journal = {Combinatorics, Probability and Computing},
  author  = {Ehard, Stefan and Glock, Stefan and Joos, Felix},
  year    = {2020},
  pages   = {868--885}
}

@article{CL06,
  author     = {Chung, Fan and Lu, Linyuan},
  title      = {Concentration inequalities and martingale inequalities: a
                survey},
  journal    = {Internet Math.},
  fjournal   = {Internet Mathematics},
  volume     = {3},
  year       = {2006},
  number     = {1},
  pages      = {79--127},
  issn       = {1542-7951,1944-9488},
  mrclass    = {60E15 (60G42)},
  mrnumber   = {2283885},
  mrreviewer = {Dominique\ L\'epingle},
  url        = {http://projecteuclid.org/euclid.im/1175266369}
}

@misc{MO,
  author       = {Sagdeev, Arsenii},
  title        = {Does any set of dominoes tile some common figure?},
  howpublished = {MathOverflow},
  year         = {2022},
  url          = {https://mathoverflow.net/questions/425311/does-any-set-of-dominoes-tile-some-common-figure},
  note         = {Posted 2022-06-29. \url{https://mathoverflow.net/q/425311}, accessed 2026-11-12}
}

@book{LP17,
  author    = {Levin, David A. and Peres, Yuval},
  title     = {{Markov} Chains and Mixing Times},
  edition   = {2},
  publisher = {American Mathematical Society},
  address   = {Providence, RI},
  year      = {2017},
  note      = {With contributions by Elizabeth L. Wilmer},
  doi       = {10.1090/mbk/107}
}

\end{document}